\documentclass[11pt]{amsart} 
\usepackage{amsfonts,amssymb,amsmath}
\usepackage{latexsym}

\numberwithin{equation}{section}

\newtheorem{thm}{Theorem}[section]
\newtheorem{prop}[thm]{Proposition}
\newtheorem{cor}[thm]{Corollary}

\newtheorem{lemma}[thm]{Lemma}
\newtheorem{rmk}[thm]{Remark}

\def\c{{\mathbb C}}
 
\def\r{{\mathbb R}}
\def\Q{{\mathbb Q}}
\def\q{{{\mathbb Q}_p}}
\def\l{{\mathbb L}}
\def\k{{\mathbb K}}

\def\O{{\rm{Ord}}}

\begin{document}

\title[On abstract homomorphisms of algebraic groups]
{On abstract homomorphisms of algebraic groups}

\author{Pralay Chatterjee}

\address{The Institute of Mathematical Sciences, C.I.T. Campus, 
Taramani, Chennai- 600113, India}

\email{pralay@imsc.res.in}

%\date{\today}

\baselineskip=16pt

\begin{abstract}
In this paper we study abstract group homomorphisms between
the groups of rational points of 
linear algebraic groups which are not necessarily reductive. 
One of our main goal is to obtain results on homomorphisms 
from the groups of rational points of
linear algebraic groups defined over certain specific fields to the groups of rational points of
linear algebraic groups over number fields and non-archimedean local fields of characteristic zero; in this set-up
we deal with the unexplored topic of abstract  homomorphisms from the groups of rational points of anisotropic groups over 
non-archimedean local fields of characteristic zero.
We also obtain results on 
abstract homomorphisms from
unipotent and solvable groups, and 
prove general results on the structures of
abstract homomorphisms using the celebrated result of Borel and Tits in this area and a well-known theorem due to Tits on the structure of the groups of rational points of isotropic semisimple groups.  
\end{abstract}
\maketitle

%\tableofcontents

\section{Introduction}
Let $G_1$ and $G_2$ be linear algebraic groups defined over fields $\l_1$ and $\l_2$, respectively, and let $G_1(\l_1)$
and $ G_2(\l_2)$ be the corresponding groups of rational points. 
In this paper we are concerned with the structures and images of
abstract group-homomorphisms $\phi : G_1(\l_1) \to G_2(\l_2)$, 
when the algebraic groups $G_1, G_2$ are not necessarily reductive.
The endeavour to understand abstract homomorphisms  as above originated in the works of
O. Schreier and B. L. van der Waerden \cite {ScV} dating back to
$1928$ and culminated in the milestone results due to A. Borel and J. Tits in $1973$,
which describes their structures under the assumptions that $G_1, G_2$ 
are absolutely simple, $G_1$ is $\l_1$-isotropic and that $\phi(G_1(\l_1))$ is Zariski-dense in
$G_2$ (see Theorem \ref{borel-tits} or \cite[Theorem A]{BT}). Since the work of Borel and Tits many important new results and generalizations have emerged; see  \cite{C}, \cite{S}, \cite{W1} for some of the advances in the recent past,
\cite{W2} for an overview and  \cite{J-W-W} for a collection of problems. 
Although the emphasis in the area seems to have been around the important case  when both $G_1$ and $G_2$ are semisimple, in \cite{BT} itself the non-semisimple situations were discussed, and it was pointed out  
through a general construction in \cite[\S 8.18, (b)]{BT} that $G_2$ may acquire a non-trivial
unipotent radical even if $G_1$ is semisimple and $\phi(G_1(\l_1))$ is Zariski-dense in $G_2$. Subsequently, a conjecture was formulated in  \cite[\S 8.19]{BT} regarding the structure of such homomorphisms for general $G_2$ under the same Zariski density condition. A proof of the conjecture is sketched by Tits (see \cite{T2}) in the special case when $\l_1=\l_2=\r$, and later it was confirmed in [L-R] with the additional hypothesis that the characteristic of $\l_1$ is zero, $G_1$ is $\l_1$-split, absolutely simple, simply connected and that the unipotent radical of $G_2$ is commutative. Very recently, along the same direction, significant results are obtained in \cite[Theorem 6.7]{Ra} considering abstract homomorphisms from the elementary subgroups of 
the groups of $R$-points of universal Chevalley-Demazure group schemes where $R$ is a
commutative ring with $1$.

Despite many contributions as above not much literature is available on the abstract homomorphisms from the groups of rational points of unipotent algebraic groups, solvable algebraic groups, anisotropic semisimple groups, or general algebraic groups which are not necessarily reductive.
While some of these questions in their full generality are very hard, one of
our main goals in this paper is to consider them in the specific set-up when $\l_2$ is either finite extensions of $\Q$ or $\q$ and $\l_1$ is either $\r,$ or an algebraically closed field, or $\Q_l$ ($l \neq p$), and describe the images $\phi(G_1 (\l_1))$ using our results in \cite{Ch1} on divisibility, orders and exponentiality in $p$-adic algebraic groups; see Theorem \ref{thm-solv}, Theorem \ref{real}, Corollary \ref{real-cor},
Theorem \ref{alg-closed}, Theorem \ref{isotropic}, Theorem \ref{sl(D)} and Theorem \ref{final-l-adic}. 
Further, the Theorem \ref{sl(D)} deals with abstract homomorphisms involving anisotropic groups which were unexplored before in the above context.
On the other hand, we establish results on abstract homomorphisms between the groups of rational points of unipotent groups; see Theorem \ref{main-unip}, Corollary \ref{main-unip-cor} and
Theorem \ref{low-card}. We also obtain general results on the structure of abstract homomorphisms appealing to the above mentioned landmark Theorem \ref{borel-tits} of A. Borel and J. Tits and the consequence Theorem \ref{kneser-tits} (4) of a fundamental result \cite[Theorem A]{T1} 
due to J. Tits; see Proposition \ref{borel-tits-app}, Corollary \ref{borel-tits-app-cor-1},
Corollary \ref{borel-tits-app-cor-2}, Corollary \ref{borel-tits-app-cor-3}, Corollary \ref{borel-tits-app-cor-4} and Theorem \ref{real-image}.
 
We now describe some of our main results.
We refer to \S 2 for the explanations on notation and terminologies used hereafter.
Henceforth, in this paper, we fix a prime number $p$ and $\k$ will denote a number field or a finite a extension of $\q$. 

Our first main result is on abstract homomorphisms $\phi : U(\l_1) \to \widetilde{U}(\l_2)$, where
$U$ and $\widetilde{U}$ are unipotent groups over fields $\l_1$ and $\l_2$, respectively,  of characteristic zero.  
We note that if $U,  \widetilde{U}$ are non-trivial unipotent groups as above then 
it is always possible to construct non-trivial abstract homomorphisms from $U(\l_1)$ to $\widetilde{U}(\l_2)$ where the image is abelian; this
can be easily done as $U(\l_1)/ [U(\l_1), U(\l_1)] $ is isomorphic to the additive group of  
a non-zero finite dimensional vector space over $\l_1$. 
Moreover, as explained in Remark \ref{unip-rmk}, if there is a ring homomorphism
from $\l_1$ to a finite dimensional algebra over $\l_2$ (or, equivalently, to a matrix algebra over $\l_2$) 
then the group of rational points of any unipotent group  
over $\l_1$ can be embedded in the group of $\l_2$-rational points of another unipotent group  over $\l_2$. 
The following result handles the opposite case when $\l_1$  never embeds in any
finite dimensional algebra over $\l_2$, and establishes that, in this case, with $U,  \widetilde{U}$ as above,
the image of any abstract homomorphism from $U(\l_1)$ to $\widetilde{U}(\l_2)$ has to be abelian, 
forcing the first typical situation mentioned above. 

\begin{thm}\label{main-unip} 
Let $\l_1,\l_2$ be fields with ${\rm Char} \l_1= {\rm Char}\l_2=0$
such that, for all $n$, the only ring-homomorphism, from $\l_1$ to
the ring ${\rm M}_n ( \l_2)$ of $n \times n$ matrices over $\l_2$,  is the trivial zero homomorphism.
Let $U$ and $\widetilde{U}$ be unipotent groups over the fields $\l_1$ and $ \l_2$, respectively.
Then for any abstract homomorphism
$\phi : U (\l_1) \to \widetilde{U}(\l_2)$ the image $\phi ( U (\l_1))$ is abelian.
\end{thm}

We show in Lemma \ref{divalg} that the field of real numbers $\r$, any algebraically closed field $\l$
with ${\rm Char}\l=0$ and the fields $\Q_l$, for primes $l\neq p$,
can never be embedded in a matrix algebra over $\q$, satisfying the hypothesis of Theorem \ref{main-unip}. Thus the following corollary follows immediately from Theorem \ref{main-unip}, Theorem \ref{unipotent} and Lemma \ref{divalg} (and hence the proof is omitted).

\begin{cor}\label{main-unip-cor} 
Let $G$ be a linear algebraic group defined over $\k$.
Let $\l$ be either an algebraically closed field with ${\rm Char}\l=0$ or the field of 
real numbers $\r$ or a finite extension of $\Q_l$ for some prime $l \neq p$, and let
$U$ be a unipotent group over $\l$. Then for any abstract homomorphism
$\phi: U (\l) \to G (\k)$ the subgroup ${\rm Zcl} \, \phi (U (\l))$ is abelian and unipotent.
\end{cor} 

In view of the above considerations on abstract homomorphisms involving unipotent groups it is now
appropriate to mention the Corollary \ref{unipotent-Q} where we prove that
if $\l_1,\l_2$ are number fields, $U$ is a unipotent algebraic group over $\l_1$ and $G$
is an algebraic group over $\l_2$ then any abstract homomorphism from $U (\l_1)$ to $G(\l_2)$ is essentially an ``algebraic'' one.

In the same spirit as of Theorem \ref{main-unip}, our next result deals with abstract homomorphisms from the group of rational points of unipotent groups to abstract groups with cardinality strictly less than that of the underlying
field over which the unipotent group is defined.

\begin{thm}\label{low-card} 
Let $H$ be an abstract group and let $\l$ be field with  ${\rm Char}\l=0$.
Further assume that cardinality of  $\l$ is strictly greater than that of 
$H$. Then for any unipotent group $U$, defined over $\l$, and for any abstract homomorphism
$\phi : U (\l) \to H$ the image $\phi ( U (\l))$ is abelian.
\end{thm}

Let $\l$ be a field with ${\rm Char} \l=0$ such that the
cardinality of $\l$ does not exceed that of $\c$, and $H$ be an algebraic group over $\l$.
We next see how easily one constructs abstract embeddings of $H (\l)$ into the groups of rational points of
certain algebraic groups, associated to $H$, which are defined over $\c$ and $\r$. 
Fixing an embedding $\alpha : \l \to \c$ we produce the algebraic group $ ^{\alpha}\! H$  over $\alpha (\l)$, 
by performing on $H$ a base change from $\l$ to $\alpha (\l)$, and assign the natural Zariski-dense embedding
$\alpha^0 \!: H (\l) \to {^{\alpha}\! H} (\c)$ induced by $\alpha$; see \S 2 for the definition of $\alpha^0$.
Now the Weil restriction ${\mathcal R}_{\c/\r} ( ^{\alpha}\! H )$ of  $^{\alpha}\! H$, which identifies
${\mathcal R}_{\c/\r} ( ^{\alpha}\! H ) (\r)$ and $^{\alpha}\! H (\c)$, in turn yields another Zariski-dense embedding of $H(\l)$ in ${\mathcal R}_{\c/\r} ( ^{\alpha}\! H ) (\r)$. 
In the following Theorem \ref{real-image} we show that, if $H$ is further assumed to be absolutely simple and $G_1, G_2$ are arbitrary
algebraic groups over $\c$ and $\r$, respectively, then arbitrary non-trivial abstract homomorphisms from $H ^+ (\l)$ to $G_1 (\c)$ or $G_2 (\r)$ are closely related to the above construction. See \S 2 for the relevant notation used in the next result.

\begin{thm}\label{real-image}
Let $\l$ be a field with $ {\rm Char} \l=0$ and let $H$ be an algebraic group over $\l$ which is $\l$-isotropic and absolutely simple.
Let $H^+ (\l)$ be the group generated by the unipotent elements of $H (\l)$. 
\begin{enumerate}
\item Let $G_1$ be an algebraic group over $\c$. Then for any non-trivial 
abstract homomorphism $\phi : H^+ (\l) \to G_1(\c)$ there exists an embedding
$\alpha : \l \to \c$ and a surjective algebraic group homomorphism 
$\mu : {\rm Zcl}\, \phi (H^+ (\l)) \to {^{\alpha}\! H}/ Z(^{\alpha}\! H)$  (over $\c$) such that
$\mu \circ \phi = \delta \circ \alpha^0$ where $\delta : {^{\alpha}\! H}
\to {^{\alpha}\! H} / Z(^{\alpha}\! H)$ is the quotient homomorphism.
In particular, we have ${\dim} (G_1/ R(G_1)) \geq  {\dim} H$, where $R(G_1)$ is the solvable radical of $G_1$.

\item Let $G_2$ be an algebraic group over $\r$. Moreover, if $\l$ does not embed in $\r$ then for any non-trivial 
abstract homomorphism $\phi : H^+ (\l) \to G_2(\r)$ there exists an embedding
$\alpha : \l \to \c$ and a surjective algebraic group homomorphism
$\mu : {\rm Zcl}\, \phi (H^+ (\l)) \to {\mathcal R}_{\c / \r} (^{\alpha}\! H/ Z(^{\alpha}\! H))$, over $\r$,
such that $\pi \circ \mu \circ \phi =  \delta \circ \alpha^0$ where
$\delta : {^{\alpha}\! H} \to ^{\alpha}\! H/ Z(^{\alpha}\! H)$
is the quotient homomorphism and
$\pi : {\mathcal R}_{\c / \r} (^{\alpha}\! H/ Z(^{\alpha}\! H))  \to {^{\alpha}\!H} / Z(^{\alpha}\! H)$
is the natural projection over $\c$ which induces an isomorphism
$\pi :{\mathcal R}_{\c / \r} (^{\alpha}\! H/ Z(^{\alpha}\! H)) (\r) \to {^{\alpha}\!H} / Z(^{\alpha}\! H) (\c)$
of abstract groups. In particular, we have 
${\dim} (G_2/ R(G_2)) \geq 2 {\dim} H$, where $R(G_2)$ is the solvable radical of $G_2$.
\end{enumerate}
\end{thm}

Recall that the field $\q$ embeds in $\c$ but not in $\r$, and thus we may apply the above result by specifying $\l = \q$.

The proof of Theorem \ref{real-image} is obtained by applying a very general result formulated 
in Proposition \ref{borel-tits-app}. Natural consequences of Proposition \ref{borel-tits-app} are then derived in 
Corollaries \ref{borel-tits-app-cor-1}, \ref{borel-tits-app-cor-2}, \ref{borel-tits-app-cor-3} and \ref{borel-tits-app-cor-4} which are interesting in their own right. We note that with certain restrictions involving the underlying fields, Corollary 
\ref{borel-tits-app-cor-2} substantially generalises the results \cite[Theorem 1]{Che} and \cite[Corollary 2.6]{Che}; see 
Remark \ref{rmk-Che}.
The proof of Proposition \ref{borel-tits-app} follows by appealing to the Theorem \ref{borel-tits} of A. Borel and J. Tits and Theorem \ref{kneser-tits} (4) due to J. Tits.
However, the hypotheses in Proposition \ref{borel-tits-app}
and Corollaries \ref{borel-tits-app-cor-1}, \ref{borel-tits-app-cor-2}, \ref{borel-tits-app-cor-3}, \ref{borel-tits-app-cor-4} are much weaker than that in Theorem \ref{borel-tits} as they
neither demand the range to be semisimple nor do they require the image to be Zariski-dense in the range. 
Consequently, the above corollaries expand 
the scope of application of Theorem \ref{borel-tits}. Thus, given the importance of Theorem 
\ref{borel-tits}, it is very surprising that these results could not be found, to the best of our knowledge, in the existing literature.
 
The Theorem \ref{real-image} above which deals with abstract homomorphisms with range as the groups of rational points of algebraic groups over $\c$ and $\r$ naturally motivates our next goal of studying 
abstract homomorphisms with range as the groups of rational points of algebraic groups over $\k$, instead of $\c$ or $\r$, and from the following results it will be evident that the behaviour of abstract homomorphisms in this setting will be in stark contrast with that in the previous one.
The following result deals with the abstract homomorphisms from the groups of rational points of solvable algebraic groups. 

\begin{thm}\label{thm-solv} 
Let $G$ be a linear algebraic group defined over $\k$.
Let $\l$ be a field, $B$ be a (Zariski-)connected solvable group over $\l$ and $\phi : B(\l)
\to G(\k)$ be an abstract homomorphism. 
\begin{enumerate}
\item If $\l$ is algebraically closed with ${\rm Char}\l =0$ then ${\rm Zcl}\,\phi (B(\l))$ is unipotent abelian.
\item If $\l= \r$ then ${\rm Zcl}\,\phi (B(\r))$ is abelian. Moreover, if
 $B (\r)^*$ denotes the connected component of $B(\r)$ in the topology induced by $\r$
then ${\rm Zcl}\,\phi (B(\r)^*)$ is unipotent and abelian.
\item If $l$ is a prime with $l \neq p$, $\l= \Q_l$ and $B$ is $\Q_l$-split then ${\rm Zcl}\,\phi (B(\Q_l))$ is abelian.
\end{enumerate}
\end{thm}

The proof of Theorem \ref{thm-solv} is obtained using Proposition \ref{prop-solv}
along with results from \cite{BS} and Corollary \ref{main-unip-cor}.

In the above set-up we next turn our focus to abstract homomorphisms from the groups of rational points of algebraic groups which are not necessarily solvable. The following result describes the images of abstract homomorphisms from general real Lie groups.

\begin{thm}\label{real}
Let $G$ be a linear algebraic group defined over $\k$.  Let $H$ be a connected real Lie group and
$\phi : H \to G(\k)$ be an abstract homomorphism. Then ${\rm Zcl}\,\phi(H)$ is a unipotent subgroup
of $G$. Moreover, if $H'$ is an algebraic group defined over $\r$ and  if $H := H' (\r)^* $ is the
connected component of $H' (\r)$ in the topology induced by $\r$ then ${\rm Zcl}\,\phi(H)$ is a unipotent and abelian subgroup of $G$.
\end{thm}

As an immediate consequence we record the following corollary. 

\begin{cor}\label{real-cor}
Let $G$ be a linear algebraic group defined over $\k$.  Let $H$ be a connected real Lie group
and $\phi : H \to G(\k)$ be an abstract homomorphism.
If $H$ is perfect then $\phi$ is trivial. In particular, if $H$ is real semisimple then $\phi : H \to G(\k)$ is trivial.
\end{cor}

We next describe images of abstract homomorphisms from an algebraic group over an algebraically closed field.

\begin{thm}\label{alg-closed} 
Let $G$ be a linear algebraic group defined over $\k$. Let $H$ be a connected algebraic group
over an algebraically closed field $\l$ and $\phi : H \to G ( \k)$ be an abstract group homomorphism.
\begin{enumerate} 
\item If ${\rm Char}\l =0$ then ${\rm Zcl}\,\phi(H)$ is a unipotent abelian subgroup of $G$.

\item If $H$ is semisimple and $\l$ is of any characteristic then $\phi$ is trivial.
\end{enumerate}
\end{thm}

We note that Theorem \ref{alg-closed} generalises \cite[Theorem A and Corollary A]{KM}.

Let now $l$ be a prime with $l \neq p$. The following results are devoted to studying the
images of abstract homomorphisms from the groups of rational points of algebraic groups over $\Q_l$
to the groups of rational points of algebraic groups over $\k$.
Th following theorem deals with abstract  homomorphisms from the groups of rational points of $\Q_l$-isotropic 
semisimple algebraic groups, for a prime $l \neq p$.

\begin{thm}\label{isotropic} 
Let $G$ be a linear algebraic group defined over $\k$. 
Let $l$ be a prime number with $l \neq p$. Let $H$ be a semisimple group over $\Q_l$ which does not have any $\Q_l$-anisotropic factor and let $\phi : H(\Q_l) \to G (\k)$ be an abstract group homomorphism.
Then $\phi ( H^+ (\Q_l)) = e$.
\begin{enumerate}
\item If $H$ is further assumed to be simply connected then $\phi ( H(\Q_l)) = e$.
\item  Moreover, $\phi ( H(\Q_l))$ is a finite abelian subgroup of $G (\k)$.
\end{enumerate}
\end{thm}

If $\k$ is a finite extension of $\q$, as both $H(\Q_l)$ and $ G (\k)$
are locally compact groups, we may consider continuous
homomorphism $\phi : H(\Q_l) \to G (\k)$ .
In \cite[Proposition 2.6.1 (ii)]{M} it is proved that if $l \neq p$ and $H$ is a (connected) simply connected semisimple group over $\Q_l$ without any $\Q_l$-
anisotropic factor then any continuous
homomorphism $\phi : H(\Q_l) \to G (\k)$ is trivial. Thus, as 
continuity of $\phi$ is not assumed, Theorem \ref{isotropic} (1) can be regarded as a  generalization of  \cite[Proposition 2.6.1 (ii)]{M}.

Abstract homomorphisms
from the groups of rational points of anisotropic groups are much less understood as, in the absence
of non-trivial unipotent elements, they defy 
techniques applied to the isotropic ones.  
Continuing in the same set-up as above we next consider the almost unexplored topic of abstract homomorphisms from the groups of rational points of anisotropic groups over $\Q_l$, for a prime 
$l \neq p$. 
As the group of rational points of an absolutely simple, simply connected and $\Q_l$-anisotropic group is isomorphic to ${\rm SL}_1 (D)$ for some central division algebra $D$ over a finite extension of $\Q_l$,  
the case of ${\rm SL}_1 (D)$ is the 
essential one that needs to be handled. 
Invoking results from \cite{Ch1} we accomplish this in the following theorem.
  
\begin{thm}\label{sl(D)} 
Let $G$ be a linear algebraic group defined over $\k$. 
Let $l$ be a prime number with $l \neq p$. Let $D$ be a central division algebra over a finite extension $\l$ of $\Q_l$. Let $r$ be the degree of the residue field of $\l$ over the residue field of $\q$ and let $d= \sqrt{{\rm dim}_\l D}$.
Let $\phi : {\rm SL}_1 (D) \to G ( \k)$ be an abstract group homomorphism. Then the image 
$\phi ( {\rm SL}_1 (D))$ is a finite solvable group with 
cardinality dividing an integer of the form $(1 + l^r + \cdots + (l^r)^{d-1}) l^m$, for some integer $m$.
\end{thm}

Our final result, which deals with
abstract homomorphisms from the groups of rational points of general algebraic groups over $\Q_l$, for primes $l \neq p$, to the groups of rational points of algebraic groups over $\k$, is proved 
combining Theorem  \ref{isotropic}, Theorem \ref{sl(D)} and Corollary \ref{main-unip-cor}..

\begin{thm}\label{final-l-adic}
Let $G$ be a linear algebraic group defined over $\k$. 
Let $l$ be a prime number with $l \neq p$. Let $H$ be an connected algebraic group over $\Q_l$ 
and $\phi : H(\Q_l) \to G (\k)$ be an abstract group homomorphism.
\begin{enumerate}
\item If $H$ is semisimple then $\phi ( H(\Q_l))$ is a finite solvable subgroup of $G (\k)$. 
\item If $H$ is any algebraic group then ${\rm Zcl}\, \phi ( H(\Q_l))$ is a solvable algebraic group. Moreover, if $H/ R_u (H)$ is semisimple then ${\rm Zcl}\,\phi ( H(\Q_l))$ is a semidirect product of a finite solvable group and an abelian unipotent radical.
\end{enumerate}
\end{thm}

The paper is organised as follows. In \S 2 we fix some standard notations and recall some well-known results.
In \S 3 we derive some preliminary results, Proposition \ref{p-divisible} and Corollary \ref{unipotent-Q},
mostly using results from \cite{Ch1}. In \S 4 we deal with homomorphisms from unipotent and solvable groups, and  Theorem \ref{main-unip}, Lemma \ref{divalg}, Corollary \ref{main-unip-cor},
Theorem \ref{low-card}, Proposition \ref{prop-solv} and Theorem \ref{thm-solv} are proved. 
In \S 5 we prove Proposition \ref{borel-tits-app}, Corollary \ref{borel-tits-app-cor-1}, 
Corollary \ref{borel-tits-app-cor-2}, Corollary \ref{borel-tits-app-cor-3}, Corollary \ref{borel-tits-app-cor-4}, Theorem \ref{real-image}, Theorem \ref{real}, Theorem \ref{alg-closed}, 
Theorem \ref{isotropic}, Theorem \ref{sl(D)} and Theorem \ref{final-l-adic}.

\section{Notation and background} 
In this section we fix the notations, definitions, recall 
some standard facts and collect some well-known results, which will be used throughout
this paper; a few specialized ones are  mentioned as and when they are required later.
The reader is referred to \cite{Sp} for generalities in the 
theory of algebraic groups.

In this paper, as mentioned in \S 1, we fix a prime number $p$ and $\k$ will always denote a number field or a finite a extension of $\q$.

The {\it cardinality} of a set $S$ is denoted by $| S|$.
For a field $\l$ we always fix an {\it algebraic closure}
which is denoted by $\overline\l$.  We denote the {\it characteristic} of
the field $\l$ by ${\rm Char}\l$. For  a pair of abelian groups $A_1, A_2$ the {\it set of homomorphisms}
from $A_1$ to $A_2$, which also acquires a natural abelian group structure,
is denoted by ${\rm Hom}_+  (A_1, A_2)$.

The {\it center} of a group $G$ is denoted by $Z(G)$.
For a group $G$ and a subset $\Lambda \subset G$ we denote the
subgroup of elements commuting with all the elements of $\Lambda$, by $Z_G (\Lambda)$. 
For an algebraic group (or a real Lie group) $E$ 
let $L(E)$ denote its {\it Lie algebra}. If $G$ is an algebraic group (or a Lie group)
and $\Lambda \subset L(G)$ then
we denote the subgroup of elements of $G$ which fix all the elements of $\Lambda$ 
via the adjoint representation, by $Z_G (\Lambda)$. 
If $V$ is an algebraic variety then the {\it Zariski closure} of a subset $W \subset V$ in $V$
is denoted by ${\rm Zcl}\,W$.  The {\it Zariski connected component of identity} in
an algebraic group $G$ is denoted by $G^0$. If $H$ is a real Lie group the
{\it connected component of identity} of $H$, in the real topology,
is denoted by $H^*$ (this non-standard notation is introduced to avoid confusion 
with the Zariski connected component). 
For an algebraic group $G$ the {\it solvable radical} 
and the {\it unipotent radical} are denoted by $R(G)$ and $R_u (G)$, respectively. If $G$ is a real Lie group
we continue to denote the {\it solvable radical} by $R(G)$.

If $G$ is an algebraic group and $x \in G$ then in the Jordan decomposition of $x \in G$ (resp. $X \in L(G)$)
the {\it semisimple} and the {\it unipotent} (resp. nilpotent) parts
of $x$ (resp. X) are denoted by $x_s$ and $x_u$ (resp. $X_s$ and $X_n$) respectively.
If $G$ is an algebraic group defined over a field $\l$ 
then the Lie algebra $L(G)$ acquires a $\l$-structure compatible
with the $\l$-structure of $G$; the {\it group $\l$-rational points} in $G$
and the {\it $\l$-Lie subalgebra of $\l$-rational points} in $L(G)$ are denoted by $G (\l)$ and
$L(G)(\l)$, respectively.
For an algebraic group $G$, defined over $\l$, the {\it variety of unipotent elements} 
of $G$ and the {\it variety of nilpotent elements} of $L(G)$ are both defined over  $\l$ and
we denote them by ${\mathcal U}_G$ and ${\mathcal N}_G$,
respectively. Moreover, if  ${\rm Char}\l =0$ 
then there is a {\it exponential map}
$\exp_G : {\mathcal{N}}_G \to {\mathcal{U}}_G$ which is a $\l$-isomorphism of $\l$-varieties.
As there is no scope for confusion, in the case when $G$ is a real Lie group, the usual {\it exponential map} 
from the Lie algebra $L(G)$ to $G$ is also denoted by
$\exp_G$. For a field $\l_1$, a finite separable extension $\l_2$ of $\l_1$ and a linear algebraic group $A$ defined over $\l_2$, the {\it Weil restriction} of $A$ over $\l_1$ 
is denoted by ${\mathcal R}_{\l_2 / \l_1} (A)$.
Recall that there is a natural {\it projection homomorphism}, say, $\pi : {\mathcal R}_{\l_2/ \l_1} (A) \to A$,
defined over $\l_2$, which induces an isomorphism
$\pi : {\mathcal R}_{\l_2/ \l_1} (A) (\l_1) \to A (\l_2) $ of abstract groups.
 
Let $A$ be an algebraic group defined over $\r$.
We recall a more refined version of Jordan decomposition, called
the {\it complete  Jordan decomposition}, available
in the group $A(\r)$ and in the Lie algebra $L(A(\r))$. 
The reader is referred to \cite[Proposition 2.4]{B1} and \cite[Theorem 7.2, p. 431]{H} for details.
An element $e \in A(\r)$ (resp. $E \in L(A(\r))$) is said to be {\it compact or elliptic} if $e$ 
(resp. $\exp_{A(\r)} (E)$) lies in a compact subgroup
of $A(\r)$. Observe that a compact element is necessarily semisimple.
A semisimple element $h \in A(\r)$ (resp. $H \in L(A(\r))$) is said to be {\it hyperbolic} if 
$h \in T(\r)^*$ (resp. $H \in L(T(\r)^*)$) for some $\r$-split torus $T$ of $A$.
For a semisimple element $s \in A(\r)$ (resp. $S \in L(A(\r))$), there is a
unique pair of elements $s_e$ (resp. $S_e$) and $s_h$ (resp. $S_h$) in $A(\r)$ (resp. in $L(A(\r))$) such that
$s_e$ (resp. $S_e$) is compact, $s_h$ (resp. $S_h$) is hyperbolic and $s=s_es_h = s_hs_e$
(resp. $S= S_e + S_h$).
For $x \in A(\r)$ (resp. $X \in L(A(\r))$), the element $(x_s)_e$ (resp.$(X_s)_e$) will be denoted by
$x_e$ (resp. $X_e$). Similarly the element $(x_s)_h$ (resp.$(X_s)_h$) will be denoted by
$x_h$ (resp. $X_h$). 
The {\it positive} part of $x$ (resp. $X$), which is defined to be $x_hx_u$ (resp. $X_h + X_n$), is denoted by
$x_p$ (resp. $X_p$). An element $y \in A(\r)$ (resp $Y \in L(A(\r))$ 
is said to be {\it positive} if $y=y_h y_u$ (resp. $Y = Y_h + Y_n$). 
Thus we have the {\it complete Jordan decomposition} in $A(\r)$ (resp. $L(A(\r))$) as follows. 
For all $z \in A(\r)$ (resp. $Z \in L(A(\r))$), there is a unique mutually
commuting triplet $z_e,z_h,z_u \in A(\r)$ (resp. $Z_e,Z_h,Z_n \in L(A(\r))$)  
such that $z_e$ (resp. $Z_e$) is compact, $z_h$ (resp. $Z_h$) is hyperbolic, $z_u$ is unipotent 
(resp. $Z_n$ is nilpotent) and $z = z_ez_hz_u$ (resp. $Z = Z_e + Z_h + Z_n$).
Moreover, written differently, for all $z \in A(\r)$ (resp. $Z \in L(A(\r))$), there are commuting elements
$z_e, z_p \in A(\r)$ (resp. $Z_e,Z_p \in L(A(\r))$) such that $z_e$ (resp. $Z_e$) is compact, $z_p$ (resp. $Z_p$)
is positive and that 
$z = z_e z_p$ (resp. $Z = Z_e + Z_p$). We will denote the set of positive elements in $A(\r)$ (resp. $L(A(\r))$)
by ${\mathcal P}_{A(\r)}$ (resp. ${\mathcal P}_{L(A(\r))}$).  We need the following lemma. 

\begin{lemma}\label{bijection-positive}{\rm (cf. \cite[Lemma 3.2]{Ch2})}
Let $G$ be a linear algebraic group defined over $\r$.
Then the exponential map $\exp_{G(\r)} : {\mathcal P}_{L(G(\r))} \to {\mathcal P}_{G(\r)}$ is a bijection.
In particular, if $g \in G(\r)$ and $X \in {\mathcal P}_{L(G(\r))}$ then $g \exp (X) g^{-1} = \exp (X)$
if and only if ${\rm Ad}(g) X = X$.
\end{lemma}

For a semisimple algebraic group $G$ over a field $\l$, 
let {\it $G^+ (\l)$ denote the subgroup of $G (\l)$ generated by the subgroups $R_u (P) (\l)$}
where $P$ varies among all the $\l$-parabolic subgroups of $G$. Recall that  $G^+ (\l)$
is a normal subgroup of $G(\l)$, and moreover, $G^+ (\l)$ coincides with the group generated by the
unipotent elements of $G(\l)$ when $\l$ is perfect.
Let $\l_1, \l_2$ be fields and let $\alpha : \l_1 \to \l_2$ be a field homomorphism.
Let $A$ be an algebraic group defined over $\l_1$, and let {\it $^\alpha\! A$ be the algebraic group defined over 
the field $\alpha(\l_1)$
obtained from $A$ by changing the base to $\alpha (\l_1)$ from $\l_1$}. We denote {\it the abstract group
homomorphism from $A (\l_1)$ to $^\alpha\! A (\l_2)$, induced from $\alpha : \l_1 \to \l_2$, by
$\alpha^0 : A (\l_1) \to {^\alpha\! A (\l_2)}$}.
We now recall the celebrated result of A. Borel and J. Tits on abstract homomorphism of algebraic groups.

\begin{thm}\label{borel-tits}{\rm (Borel-Tits, \cite[Theorem A]{BT})}
Let $\l_1, \l_2$ be infinite fields. Let $G_1$ and $G_2$ be algebraic groups defined over
$\l_1$ and $\l_2$, respectively. Let $H \subset G_1 ( \l_1)$ be a subgroup containing $ G_1^+ (\l_1)$
and let $\phi : H \to G_2 (\l_2)$ be an abstract homomorphism of groups. Assume that
\begin{enumerate}
\item $G_1$ is $\l_1$-isotropic.
\item $G_1$ and $G_2$ are absolutely simple.
\item either $G_1$ is simply connected or $G_2$ is adjoint.
\item $G_2 \neq 1$ and ${\rm Zcl}\,\phi (G^+_1(\l_1)) = G_2$.
\end{enumerate}
Then there exists a unique field homomorphism $\alpha : \l_1 \to \l_2$,
an isogeny $\beta :  {^\alpha \! G_1} \to G_2$ of $\l_2$-algebraic groups and a
group homomorphism $\tau : H \to Z (G_2) (\l_2)$ such that 
$$
\phi (g) = \tau (g) \beta (\alpha^0 (g)), \,\, for \,\, all \,\, g \in G_1( \l_1).
$$
\end{thm}

We now collect certain well-known facts regarding the subgroup $G^+ (\l)$ in the case when $G$ is semisimple. 

\begin{thm}\label{kneser-tits} 
Let $\l$ be a field and let $G$ be a semisimple algebraic group over $\l$. 
Assume further that $G$ has no $\l$-anisotropic factors.
\begin{enumerate}
\item {\rm (cf. \cite[p. 406]{PR})} If ${\rm Char} \l =0$ and $\pi: \widetilde{G} \to G$ is a (algebraic) covering homomorphism over $\l$ then $\pi (\widetilde{G}^+ (\l)) =  G^+ (\l)$. 
\item {\rm (cf. \cite[Theorem 2.3.1, p. 52]{M})} $G (\l) / G^+ (\l)$ is a finite abelian group.
\item {\rm (Platonov, cf. \cite[Theorem 2.3.1, p. 52]{M})} If $\l$ is a non-archimedean local field and $G$ is simply connected then $G (\l) = G^+ (\l)$.
\item {\rm (Tits \cite{T1})} If $|\l| \geq 4 $ and $G$ is $\l$-simple then $[G^+ (\l), G^+ (\l)] = G^+ (\l)$.
\end{enumerate}
\end{thm}

We recall that (3) of Theorem \ref{kneser-tits} is regarded as Kneser-Tits conjecture which was proved by V. Platonov, and (4) is a immediate consequence of the main 
theorem of the paper \cite{T1} by J. Tits.
 
We now state the well-known Jacobson-Morozov theorem.

\begin{thm}\label{jacobson-morozov}{\rm (Jacobson, Morozov, cf. \cite[Lemma 2.8]{Wi})}
Let $\l$ be a field with ${\rm Char} \l =0$ and let $G$ be a semisimple algebraic group over $\l$. 
Let $ u \in G(\l)$ be a unipotent element. Then there is a $\l$-homomorphism $\delta : {\rm SL}_2 (\overline{\l}) \to G$
such that  
$$
\delta ( 
\begin{pmatrix}
1 & 1\\
0 & 1
\end{pmatrix}) = u.
$$
\end{thm}

The following basic result on the structure of $G(\l)$, where $\l$ is a non-discrete locally compact field
and $G$ is a $\l$-simple group, is due to C. Riehm and it will be used in the last section.

\begin{thm}\label{riehm}{\rm (Riehm \cite{Ri}, cf. \cite[Theorem 3.3, p. 114]{PR})}
Let $\l$ be a non-discrete locally compact field with
${\rm Char} \l =0$ and $G$ be a $\l$-simple algebraic group. Then every non-central
normal subgroup of $G (\l)$ is open in $G(\l)$.
\end{thm}

Finally, we need some facts associated to 
profinite groups; see \cite{R-Z}
for generalities. A topological group $G$ is called {\it profinite} 
if it is compact and totally disconnected.
Such a group is topologically isomorphic to a projective limit of
finite groups. Moreover, a profinite group is called {\it prosolvable} if it is the projective limit of finite solvable groups.
A {\it supernatural number} is a formal infinite
product $\prod p^{n(p)}$, over all primes $p$, 
where $n(p)$ is a non-negative integer
or infinity. Product, divisibility, l.c.m and g.c.d of
a set (possibly infinite) of supernatural numbers are defined in the
natural way. In particular, l.c.m of an infinite set of integers is
a supernatural number.  
We now define the {\it order}, ${\rm Ord}(G)$ of a profinite group $G$ by,
$$
{\rm Ord}(G)= 
{\rm l.c.m} \{ |G/U| \, : \,\, U \,\, {\rm is} \,\, {\rm an} \,\,
{\rm open}\,\, {\rm subgroup} \,\, {\rm of}
\,\, G \}
$$
Thus, if $G$ is a profinite group then $\O(G)$ 
is a supernatural number.
A profinite group $G$ is called a {\it pro-$p$ group} 
if $p$ is the only prime which divides $\O(G)$.
We need the following results from \cite{R-Z}.

\begin{thm}\label{profinite}{\rm (cf.  \cite[Proposition 4.2.1, p. 124]{R-Z})}
Let $G$ be a profinite group and let $N$ be a subgroup
(not necessarily closed) of $G$ of finite index. Then the index $[G : N]$ divides ${\rm Ord}(G)$.
\end{thm}

\begin{thm}\label{prosolv}{\rm (cf.  \cite[Corollary 4.2.2, p. 125]{R-Z})} 
Let $G$ be a prosolvable group and let $N$ be a normal subgroup
(not necessarily closed) of $G$ of finite index. Then $G/N$ is a finite solvable group.
\end{thm}

\section{Preliminary results}
In this section we deduce some preliminary results which will be used in the subsequent sections.
For a group $H$ let $P_n : H \to H$ denote the {\it $n$-th power map} defined by $P_n (h) = h^n, \, h \in H$.
Recall that, for an algebraic group $G$, the variety of unipotent elements of $G$ is denoted ${\mathcal U}_G$ 
and that throughout the paper $\k$ stands for a number field or a finite extension of $\q$ for a fixed prime $p$.

We require the following key results from \cite{Ch1}.

\begin{thm}\label{key1} {\rm (\cite[Theorem 1.5, Corollary 1.6, Theorem 6.2]{Ch1})}
Let $G$ be a linear algebraic group defined over $\k$. 
\begin{enumerate}
\item
Let $\alpha \in G (\k)$ be a semisimple element with
$\alpha \in \bigcap_{n=1}^{\infty} P_{p^n} ( G(\k))$ then
the cyclic subgroup generated by $\alpha$ is finite and its order is relatively prime to $p$.
\item 
$$
\bigcap_{n=1}^{\infty} P_n (G(\k)) = {\mathcal U}_G (\k).
$$

\item
Let $\phi : \Q \to G(\k)$ be an abstract group homomorphism from the additive group of rational numbers to  $G(\k)$. Then   
there is a nilpotent element $X \in L(G) (\k)$ such that 
$$
\phi (t) \, = \, \exp_G (tX), \,\, \text{for all}\, \, t \in \Q.
$$
\end{enumerate}
\end{thm}

The proofs of the above results are given for the cases $\k = \q, \Q$ in  \cite[Theorem 1.5, Corollary 1.6, Theorem 6.2]{Ch1}. 
However, as mentioned in \cite[Remark 6.3]{Ch1}, similar proofs carry over when $\k$ is a finite extension of either $\Q$ or $\q$. 

We will now record some immediate applications of the above results.

\begin{thm}\label{unipotent} 
Let $G$ be a linear algebraic group defined over $\k$. Let $H$ be an abstract group and 
$\phi : H \to G(\k)$ be an abstract group homomorphism. Let $h$ be a divisible element in $H$.
Then $\phi (h)$ is a unipotent element in  $G(\k)$.
In particular, if $H$ is an algebraic group defined over a field $\l$ with ${\rm Char}\l =0$
and $u \in H(\l)$ is a unipotent element then so is $\phi (u)$. Consequently, if $U$ is a unipotent group
over $\l$ and $\phi : U (\l) \to G (\k)$ is an abstract homomorphism then
${\rm Zcl}\, \phi (U (\l))$ is a unipotent group.
\end{thm}

\noindent{\bf Proof.}
Let $h$ be a divisible element in $H$. Then so is $\phi (h)$ in $G(\k)$. Hence, using Theorem \ref{key1},
$$
\phi (u) \in \bigcap_{n} P_n (G(\k)) = {\mathcal U}_G (\k).
$$
Thus $\phi(u) \in G(\k)$ is unipotent. 

Now to prove the second assertion, we first see that, as $u \in H(\l)$ is unipotent, $u$ lies in a
one dimensional unipotent $\l$-subgroup $U$ of $H$. 
The proof is concluded by applying the first part and noting that, as $U(\l)$ is isomorphic to the additive group
$\l$, the element $u$ is divisible in $H(\l)$.  The last statement follows immediately from the
second one.
\hfill$\Box$

\begin{rmk}
{\rm Let $G_1, G_2$ be algebraic groups over fields $\l_1$ and $\l_2$, respectively. Let
$\phi : G_1 (\l_1) \to G_2 (\l_2)$ be an abstract homomorphism and $u \in G_1 (\l_1)$ be a unipotent element.
Assuming $\l_1$ to be infinite and $G_1$ to be semisimple 
it was proved in \cite[Proposition 7.2 (i)]{BT} that $\phi (u)$ is a unipotent element in $G_2 (\l_2)$. 
However, if we replace $\l_2$ by the field $\k$, as in Theorem \ref{unipotent}, and if $\l_1$ is a field with 
${\rm Char}\l_1 =0$ then for any general algebraic group $G_1$, not necessarily semisimple, it follows that
$\phi (u)$ is a unipotent element in $G_2 (\k)$. This follows easily using
Theorem \ref{unipotent} and the fact that $u$ is a divisible element in $G_1 (\l_1)$.}
\end{rmk}

We omit the proof of the next corollary as the proof is easily obtained from Theorem \ref{key1} (2).

\begin{cor} 
Let $G$ be a linear algebraic group defined over $\k$. 
Further assume that $G$ has a Levi $\k$-subgroup which is $\k$-anisotropic.
Let $M$ be a perfect group generated by divisible elements then any abstract
homomorphism $\phi : M \to G(\k)$ is trivial. In particular, if $\l$ is field with ${\rm Char}\l =0$ and $H$ is
a semisimple algebraic group defined over $\l$
then any abstract homomorphism $\phi : H^+ ( \l) \to G(\k)$ is trivial.
\end{cor}

The following proposition will be used in the proofs of Proposition \ref{prop-solv} and Theorem \ref{sl(D)}.

\begin{prop}\label{p-divisible} 
Let $G$ be a linear algebraic group defined over $\k$. 
Let $H$ be an abstract group which is $p$-divisible and 
$\phi : H \to G(\k)$ be an abstract group homomorphism. Then there is a unipotent $\k$-subgroup $U$ of $G$ 
and a finite $p$-divisible subgroup $F$ of $G(\k)$ such that $F$ normalises $U$ and $ \phi (H) \subset F U (\k)$.
Moreover, $F$ can be chosen so that  $\pi (\phi (H) ) = F$, where
$\pi : FU \to F$ denotes the natural projection homomorphism from the $\k$-subgroup $FU$. 
\end{prop} 

\noindent{\bf Proof.}
We first observe that that there is a positive integer $r$ such that if $\Gamma \subset G (\k)$ 
is any finite subgroup then $|\Gamma|$ divides $r$. 
Note that as $\k$ is either a finite extension of $\Q$ or a finite extension of $\q$, for some prime number $p$,
using the Weil restriction functor, we may always embed $G(\k)$ in $G'(\q)$ where $G'$ is an algebraic group defined over $\q$.
Hence to prove it is enough to assume that $\k = \q$. 
Now an argument, using \cite[Theorem 1, p. 124]{Se1}, as given in the proof of \cite[Lemma 5.3]{Ch1} implies the above assertion.

We next $\k$-embed $G$ in a suitable ${\rm GL}_n ( \overline{\mathbb K })$ which is equipped with
the usual $\k$-structure so that group of $\k$-points is ${\rm GL}_n ({\mathbb K })$.
Now as every element of $H$ is $p$-divisible one has 
$$\phi (H)  \subset
\bigcap_{n=1}^{\infty} P_{p^n} ( G(\k)).$$

Using Theorem \ref{key1} it follows that if $h \in H$ then $\phi (h)_s$ is of finite order which is relatively prime
to $p$. We now appeal to the above fact to conclude that there is a  integer $l$, relatively prime to $p$,
such that the order of the element $\phi (h)_s$ divides $l$. In particular, as  matrices in 
${\rm GL}_n ( \overline{\mathbb K })$, each element $\phi (h)$ for $h \in H$ satisfy the equation :
$$
(\phi(h)^l - {\rm Id})^n = 0.
$$

In particular, all the semisimple elements in ${\rm Zcl}\, \phi (H)$ are of finite order and moreover, the order of each
semisimple elements divide the integer $l$. Let $F$ be a Levi $\k$-subgroup of ${\rm Zcl}\, \phi (H)$. 
This implies that ${\rm Zcl}\, \phi (H)= F R_u ( {\rm Zcl}\, \phi (H))$ where the product on the right hand side is a semidirect product. 
Since all semisimple elements of $F$ are of finite order
it follows that $F$ is a finite group.
We set $U :=R_u ( {\rm Zcl}\, \phi (H))$. Let $\pi : FU \to F$ be the natural projection homomorphism as in the statement of the theorem. If $A:= \pi (\phi (H) )$ then $A \subset F(\k)$. Moreover, $\phi (H) \subset A (U(\k))$.
As $A$ is finite ${\rm Zcl}\, \phi (H) = A U$. Thus $AU = FU$ and in particular $A = F$.
This completes the proof.
 \hfill$\Box$ 
  
Let $A$ be an abelian group and $\l$ be a field. Recall that, considering $\l$ as a group under the
addition, the space of group homomorphisms from $A$ to $\l$ is denoted by
${\rm Hom}_{+} ( A,\l)$. 

We will need the next lemma which is easy to prove.

\begin{lemma}\label{1-ps}
Let $A$ be an abelian group and $\l$ be a field such that ${\rm Char} \l =0$.
Let $U$ be a unipotent group over $\l$ and
$\phi: A \to U (\l) )$ be an abstract group homomorphism.
Let  $ X_1, \cdots, X_n$ be a $\l$-basis of 
Lie algebra $L ({\rm Zcl}\, \phi(A)) (\l)$.
Then there are finitely many homomorphisms $\lambda_1, \cdots, \lambda_n$ in
${\rm Hom}_{+} ( A, \l)$ such that, for all $t \in A$, 
$$
\phi(t) = \exp_U ( \sum_{i=1}^{n} \lambda_i (t) X_i).
$$  
\end{lemma}

\noindent{\bf Proof.} 
As ${\rm Zcl}\, \phi(A)$ is abelian and $U$ is unipotent the exponential map 
$\exp_U : L({\rm Zcl}\, \phi(A) ) \to {\rm Zcl}\, \phi(A)$ is a $\l$-isomorphism of the
abelian groups $L({\rm Zcl}\, \phi(A))$ and ${\rm Zcl}\, \phi(A)$. It is clear that there exist functions
$\lambda_i : A \to \l,\, i = 1, \cdots, n$ such that 
$\phi(t) = \exp_U ( \sum_{i=1}^{n} \lambda_i (t) X_i)$, for all $t \in A$. As   
$\exp_U : L({\rm Zcl}\, \phi(A) ) \to {\rm Zcl}\, \phi(A)$ is a $\l$-isomorphism it follows immediately that
$\lambda_i \in {\rm Hom}_{+} ( A, \l)$, for all $i= 1, \cdots, n$.
\hfill$\Box$

\begin{cor}\label{1-ps-K} 
Let $G$ be a linear algebraic group defined over $\k$. Let $A$ be an abelian divisible group.
Let $\phi: A \to G(\k)$ be an abstract group homomorphism from the group $A$
to the group $G(\k)$. Then there exists $\k$-linearly independent nilpotent elements $ X_1, \cdots, X_n$ in $L(G)(\k)$
with $[X_i, X_j]=0$, for all $i,j$, and  
finitely many homomorphisms $\lambda_1, \cdots, \lambda_n$ in
${\rm Hom}_{+} ( A,\k)$ such that, for all $t \in A$, 
$$
\phi(t) = \exp_G ( \sum_{i=1}^{n} \lambda_i (t) X_i).
$$ 
\end{cor}

\noindent{\bf Proof.} 
We first use Theorem \ref{unipotent} to see that $U:= {\rm Zcl}\, \phi(A)$ is a unipotent subgroup
of $G(\k)$ over $\k$. It is immediate that $U$ is an abelian group. We now apply Lemma \ref{1-ps}
to conclude the proof.
\hfill$\Box$

We next record an immediate application of the above result. 

\begin{cor}\label{unipotent-Q}
Let $\l_1$ and $\l_2$ be finite extensions of $\Q$, and 
let $U$ and $G$ be algebraic groups over $\l_1$ and $\l_2$, respectively. 
Assume that $U$ is unipotent. Let $\phi: U(\l_1) \to G(\l_2)$ be an abstract group homomorphism.
\begin{enumerate}
\item If $ \l_1= \l_2 = \Q$ then there is a unique algebraic $\Q$-homomorphism $\psi : U \to G$ 
such that $\phi$ is the restriction of $\psi$ to $U (\Q)$.
\item There exists a unique $\Q$-homomorphism $\psi' : {\mathcal R}_{\l_1/ \Q} (U) \to {\mathcal R}_{\l_2/ \Q} (G)$ 
such that $\phi (x) = \pi_2 \circ \psi' \circ {\pi_1}^{-1} (x), \, x \in U (\l_1)$, where
$\pi_1 : {\mathcal R}_{\l_1/ \Q} (U) \to U$ is the natural projection over $\l_1$ inducing an isomorphism of groups
$\pi_1 : {\mathcal R}_{\l_1/ \Q} (U) (\Q) \to U (\l_1)$ , and analogously, 
$\pi_2 : {\mathcal R}_{\l_1/ \Q} (G) \to G$ is the natural projection over $\l_2$ inducing an isomorphism of groups
$\pi_2 : {\mathcal R}_{\l_1/ \Q} (G) (\Q) \to G (\l_2)$.
\end{enumerate}
\end{cor}

\noindent{\bf Proof.} 
We observe that (2) follows immediately from (1) and hence we proceed with a proof for (1).
Let $U':= {\rm Zcl} \phi (U(\Q))$. Then it is clear from (2) of Theorem \ref{key1} that $U'$ is a unipotent subgroup of $G$ which is defined over $\Q$.
Thus, without loss of generality we may assume that $G$ is a unipotent group.
Let $n:= \dim U$. Recall that there is $\overline{\Q}$-basis $Y_1, \cdots , Y_n$ in $L(U) (\Q)$ for the Lie algebra $L(U)$ such that the map
$E : L(U) \to U$ defined by 
$$
E ( z_1 Y_1 + \cdots + z_n Y_n) : = \exp_U (z_1 Y_1) \cdots \exp_U (z_n Y_n), \,\, (z_1, \cdots, z_n) \in {\overline{\Q}}^n
$$
is an $\Q$-isomorphism of $\Q$-varieties $L(U)$ and $U$. 
In particular, we have 
$$U (\Q) = \{ \exp_U (t_1 Y_1) \cdots \exp_U (t_n Y_n) | \,\, (t_1, \cdots, t_n) \in {\Q}^n \,\}.
$$
Now by Corollary \ref{1-ps-K} for each $Y_i$ there exists a  nilpotent element
$X_i \in L(G) (\Q)$ such that $\phi ( \exp_U (sY_i )) = \exp_G ( s X_i)$, for all $s \in \Q$. 
We next define the map $\psi : U \to G$ by 
$$
\psi (E ( z_1 Y_1 + \cdots + z_n Y_n)) : = \exp_G (z_1 X_1) \cdots \exp_G (z_n X_n), \,\, (z_1, \cdots, z_n) \in {\overline{\Q}}^n.
$$
As $E : L(U) \to U$ is a $\Q$-isomorphism of $\Q$-varieties it follows that $\psi : U \to G$ is a $\Q$-morphism of varieties defined over $\Q$. Moreover, it is clear that
$\phi$ is the restriction of $\psi$ to the subgroup $ U (\Q)$. Finally it remains to prove that $\psi$ is a group homomorphism.
This follows immediately taking into account
that $\phi : U (\Q) \to G (\Q)$ is a group homomorphism and the fact that
$U (\Q)$ is Zariski-dense in $U$. The uniqueness of $\psi$ also follows from 
the Zariski-density of $U (\Q)$ in $U$.
\hfill$\Box$

In \cite[\S 6, p.5]{St} it is proved that an abstract homomorphism $\phi : {\rm SL}_n (\Q) \to {\rm GL}_m ({\c})$ is polynomial,
i.e., it is the restriction of an algebraic homomorphism from ${\rm SL}_n (\c)$ to 
${\rm GL}_m ({\c})$. In view of this fact and
Corollary \ref{unipotent-Q} we have the following result.

\begin{cor} Let $G$ and $H$ be algebraic group over $\Q$. Assume that $H/R_u (H)$ is isomorphic (as $\Q$-algebraic groups) 
to ${\rm SL}_n$ which is equipped with the usual $\Q$-structure. Then for any abstract homomorphism $\phi : H (\Q) \to G (\Q)$
there is a unique algebraic $\Q$-homomorphism $\psi : H \to G$ 
such that $\phi$ is the restriction of $\psi$ to $H (\Q)$. 
\end{cor}

\section{Homomorphisms from unipotent and solvable groups}
In this section we prove Theorem \ref{main-unip} on the abstract homomorphisms between unipotent groups over general fields of characteristic zero. In view of Theorem \ref{main-unip} we then derive Lemma \ref{divalg} 
to obtain Corollary \ref{main-unip-cor}. Theorem \ref{low-card} is deduced independently in this section.   
We prove Proposition \ref{prop-solv} 
to deal with abstract homomorphism from solvable algebraic groups over certain specific fields, which is 
then applied to establish Theorem \ref{thm-solv}.
Recall that $\k$ stands for a number field or 
a finite extension of $\q$ for a fixed prime $p$.

We begin with some simple algebraic facts.
Let $\l_1,\l_2$ be fields. Recall that 
${\rm Hom}_{+} ( \l_1,\l_2)$ stands for the space of abstract 
homomorphisms from $\l_1$ to $\l_2$ when $\l_1, \l_2$ are considered
as abelian groups with respect to the additions. 
We define the obvious  actions of the fields $\l_1, \l_2$ on the abelian group ${\rm Hom}_{+} ( \l_1,\l_2)$ which make the latter vector spaces over both $\l_1, \l_2$.
For $x \in \l_1, y \in \l_2$ and $\lambda \in {\rm Hom}_{+} ( \l_1,\l_2)$ we define
$x \cdot \lambda$ and $y \cdot \lambda$ by $ x \cdot \lambda (z) :=  \lambda (xz)$
and $y \cdot \lambda (z) := y \lambda (z)$ for all  $ z \in \l_1$.
Note that for the sake of simplicity we used the same notations to denote both the
actions. Clearly the actions commute, that is, $ x \cdot (y \cdot \lambda) = y \cdot (x \cdot \lambda) $ for all $x \in \l_1, y \in \l_2$ and $\lambda \in {\rm Hom}_{+} ( \l_1,\l_2)$. In particular, this yields a 
a natural ring homomorphism
$\Phi : \l_1 \to {\rm End}_{\l_2} ( {\rm Hom}_+ ( \l_1,\l_2))$ given by
$\Phi (s) (\lambda) := s \cdot \lambda, \, s \in \l_1, \lambda \in {\rm Hom}_+ ( \l_1,\l_2) $.

\begin{lemma}\label{additive-homom}
Let $\l_1,\l_2$ be a pair of fields such that, for all $n$, the only ring homomorphism from $\l_1$ to
the ring ${\rm M}_n ( \l_2)$ of $n \times n$ matrices over $\l_2$ is the trivial zero homomorphism.
Let $\delta, \lambda_i, \mu_i; \,\, 1 \leq i \leq n$ be elements in ${\rm Hom}_{+} (  \l_1, \l_2)$  satisfying the relation
$$
\delta (st) = \sum _{i=1}^n \lambda_i (s) \mu_i (t) 
$$
for all $s,t \in \l$. Then $\delta =0$.
\end{lemma}

\noindent{\bf Proof.} We consider the action of $\l_1$ on ${\rm Hom}_+ ( \l_1,\l_2)$ as defined above.
Let $W_1$ be the $\l_2$-span of the set $\{\, s \cdot \delta \, | \, s \in \l_1 \, \}$.
Then the hypothesis imposed on $\delta$ implies that
$W_1$ is contained in the $\l_2$-span of finitely many elements $\mu_1, \cdots ,\mu_n$. In particular,
${\rm dim}_{\l_2} W_1 < \infty$. Further $W_1$ is a $\l_1$-invariant subspace of the $\l_2$ vector space
${\rm Hom}_{+} ( \l_1 , \l_2)$.
Since ${\rm dim}_{\l_2} W_1 < \infty$ the ring homomorphism
$\Phi : \l_1 \to {\rm End}_{\l_2} ( {\rm Hom}_{+} ( \l_1,\l_2))$ induces a ring homomorphism $\Phi' :  \l_1 \to 
{\rm End}_{\l_2}  (W_1)$, by restricting
endomorphism $\Phi (\cdot)$ to $W_1$. Let $l := {\rm dim}_{\l_2} W_1$.
As a ring ${\rm End}_{\l_2}  (W_1)$ is isomorphic to ${\rm M}_l (\l_2)$ and hence
$\Phi' (  \l_1) =0$. This implies that $s \cdot \delta =0$, for all $s \in \l$. In particular,
$\delta =0$.
\hfill$\Box$

We now give a proof of Theorem \ref{main-unip}.

\noindent
{\bf Proof of Theorem \ref{main-unip}.}
We will prove the theorem by induction on the dimension of $U$. If ${\rm dim}\, U=1$ then the group $U (\l_1)$ is abelian and the statement follows trivially.
Assuming that the statement holds for all unipotent groups $U'$ (defined over $\l_1$) with  ${\rm dim}\,U' < {\rm dim}\,U$
we will prove it for $U$.
As $U$ is $\l_1$-split there are normal $\l_1$-closed subgroups $U_i, \, 0 \leq i \leq n$ of $U$ such that
$$
U= U_0 \supset U_1 \supset \cdots \supset U_n = \{ e \}
$$
with the property that 
$U_i/U_{i+1}$ is $\l$-isomorphic to the one dimensional standard additive group ${\mathbb G}_a$ 
(defined over $\l$)
and $[U, U_i] \subset U_{i+1}$ for all $i= 0, \cdots, n-1$.
This implies that there is a element $X_0 \in L(U)(\l_1)$ such that if $E:={\rm Zcl}\, \{\exp_U (tX_0) \, | \, t \in \l_1 \}$ then 
$U(\l_1) = E(\l_1) U_1 (\l)$. Now as ${\rm dim}\,U_1 < {\rm dim }U$, by induction hypothesis, the image $\phi (U_1 (\l))$
is an abelian subgroup of $\widetilde{U} (\l_2)$. As $[U_1, U_1] (\l_1) = [U_1(\l_1), U_1(\l_1)]$ (see
\cite[Proposition 1.2]{F}) it follows that $\phi ( [U_1, U_1] (\l_1)) = \{ e \}$. 
Consequently, the abstract homomorphism $\phi : U (\l_1) \to \widetilde{U}(\l_2)$ factors through a homomorphism 
$\tilde{\phi}: U(\l_1)/ [U_1, U_1] (\l_1) \to \widetilde{U} (\l_2)$. By \cite[Corollary 1.3.1]{F} we have that
$U(\l_1)/ [U_1, U_1] (\l_1)$ is isomorphic to the group $(U/ [U_1, U_1]) (\l_1)$. 
We now divide the proof in two cases.

{\it Case 1:} We first consider the case when ${\rm dim}\, [U_1, U_1] > 0$. In this case, as ${\rm dim}\, U/ [U_1, U_1] <  {\rm dim}\, U$,
by induction hypothesis $\phi (  U (\l_1)  )= \tilde{\phi}(U(\l_1)/ [U_1, U_1] (\l_1))$ is an abelian group.

{\it Case 2:} We next consider the case when ${\rm dim}\, [U_1, U_1] = 0$. Thus $U_1$ is an abelian group.
Setting $N := E U_2$ we see that ${\rm dim}\, N = {\rm dim } U-1$. Then, by induction hypothesis, $\phi (N(\l_1)) =
\phi ( E(\l_1) U_2 (\l_1))$ is an abelian subgroup of $\widetilde{U}(\l_2)$. In particular,
$\phi ( [N(\l_1), N(\l_1)]) = \{e \}$. Further, as $[E (\l_1), U_1 (\l_1)] \subset U_2 (\l_1)
\subset N(\l_1)$, it follows that 
\begin{equation}\label{u-1}
\phi ( [ E(\l_1), [ E(\l_1), U_1 ( \l_1)] \,] )\subset \phi ([N(\l_1), N(\l_1)]) = \{e \}.
\end{equation}

Let $\widetilde{E} := {\rm Zcl}\,\phi(E(\l_1))$ and $\widetilde{U}_1 := {\rm Zcl}\,\phi(U_1 ( \l_1))$.
It is clear that $\widetilde{E},  \widetilde{U}_1$ are $\l_2$-closed abelian subgroups of $\widetilde{U}$.
Taking Zariski closure in the group $\widetilde{U}$, from the equality as in \eqref{u-1},  it follows that
\begin{equation}\label{u0}
[\widetilde{E},\,[\widetilde{E},\widetilde{U}_1 ]\,] = \{e \}.
\end{equation}
Let $v \in U_1 ( \l_1)$ be an arbitrary element. As $E(\l_1)$
and  $U_1 ( \l_1)$ are abelian, to prove that $\phi ( U (\l_1))$ is abelian it 
is enough to show that the subgroup $\phi ( E(\l_1))$ commutes with $\phi(v)$.
Let $Y \in L (U_1) (\l_1)$ be the (unique) element such that $ v = \exp_U (Y)$.

Let $p:= \dim {\rm Zcl}\,\phi(E(\l_1))$ and $q:=\dim{\rm Zcl}\,\phi(U_1 ( \l_1))$. 
Let $\widetilde{X}_1, \cdots \widetilde{X}_p$ be a $\l_2$-basis of $L (\widetilde{E}) (\l_2)$ and
$\widetilde{Y}_1, \cdots, \widetilde{Y}_q$ be a $\l_2$-basis of $L(\widetilde{U}_1)(\l_2)$.
We now apply Lemma \ref{1-ps} to see that there exist
elements $\lambda_1, \cdots, \lambda_p$ and $\mu_1, \cdots, \mu_q$ in ${\rm Hom}_+ ( \l_1,\l_2)$
such that, for all $t, s \in \l_1$,
\begin{equation}\label{u1}
\phi( \exp_U ( tX_0)) = \exp_{\widetilde{U}} ( \sum_{i=1}^p \lambda_i (t) \widetilde{X}_i),\\\
\phi( \exp_U ( sY)) =  \exp_{\widetilde{U}} (\sum_{j=1}^q \mu_j (s) \widetilde{Y}_j).
\end{equation}

Similarly, as $L(E)$ normalises $L(U_1)$, there exist $\delta_1, \cdots,   \delta_q$ in ${\rm Hom}_+ ( \l_1, \l_2)$ such that
for all $s' \in \l_2$,
\begin{equation}\label{u2}
\phi( \exp_U ( s'[X_0,Y])) = \exp_{\widetilde{U}} (\sum_{k=1}^q \delta_k (s') \widetilde{Y}_k).
\end{equation}

Now we use \eqref{u1} to see that
\begin{equation}\label{u3}
\phi( \exp_U (t X_0)) \, \phi( \exp_U (s Y)) \, \phi( \exp_U (t X_0))^{-1}
= \exp_{\widetilde{U}} ({\rm Ad}(\exp_{\widetilde{U}} ( \sum_{i=1}^p \lambda_i (t) \widetilde{X}_i))  \sum_{j=1}^q \mu_j (s) \widetilde{Y}_j).
\end{equation}

Now from \eqref{u0} we see that $[L(\widetilde{E}), \, [L(\widetilde{E}),L(\widetilde{U}_1) ]\,] = 0$. Thus we have,

\begin{equation}\label{u3-0}
[\sum_{i=1}^p \lambda_i (t) \widetilde{X}_i, \, [ \sum_{i=1}^p \lambda_i (t) \widetilde{X}_i, \sum_{j=1}^q \mu_j (s) \widetilde{Y}_j]\,]=0.
\end{equation}

Using the above equation \eqref{u3-0} we see that
\begin{equation}\label{u3-1}
\exp_{\widetilde{U}} ({\rm Ad}(\exp_{\widetilde{U}} ( \sum_{i=1}^p \lambda_i (t) \widetilde{X}_i))  \sum_{j=1}^q \mu_j (s) \widetilde{Y}_j)=\nonumber\\
\exp_{\widetilde{U}} (\sum_{j=1}^q \mu_j (s) \widetilde{Y}_j + [ \sum_{i=1}^p \lambda_i (t) \widetilde{X}_i, \, \sum_{j=1}^q \mu_j (s) \widetilde{Y}_j]).
\end{equation}

As before, since $L(\widetilde{E})$ normalises $L(\widetilde{U}_1)$, it follows that, for all $i,j$ with $1 \leq i \leq p, 1 \leq j \leq q$ 
there exist $q$ elements $a_{ijl} \in \l_2, \, 1 \leq l \leq q$ such that $[\widetilde{X}_i, \widetilde{Y}_j] = \sum_{l=1}^q a_{ijl}\widetilde{Y}_l$. 

Hence from the equalities \eqref{u3}, \eqref{u3-1} we obtain

\begin{equation}\label{u4}
\phi( \exp_U (t X_0))\, \phi( \exp_U (s Y)) \, \phi( \exp_U (t X_0))^{-1} \,
=\exp_{\widetilde{U}} ( \sum_l \{ \mu_l (s) + \sum_{i,j} a_{ijl} \lambda_i (t) \mu_j(s)\} \widetilde{Y}_l).
\end{equation}

Now a simple computation of the left hand side of \eqref{u4} yields 

\begin{gather}
\phi( \exp_U (t X_0)) \, \phi( \exp_U (s Y)) \, \phi( \exp_U (t X_0))^{-1}
=\phi ({\rm Ad} (\exp_{U} (tX_0)) sY)\nonumber\\
\label{u4-1}= \phi(\exp_U  (sY + ts [X_0,Y] + \sum_{n \geq 1} ( t^n s / n!) {\rm ad}(X_0)^n Y)).
\end{gather}

As $U_1$ is abelian the right-hand side of the above equation can be written as
\begin{equation}\label{u5}
\phi(\exp_U ( (sY + ts [X_0,Y])) \, \phi ( \exp_U ( \sum_{n \geq 1} ( t^n s / n!) {\rm ad}(X_0)^n Y)).
\end{equation}

Observe that $\sum_{n \geq 1} ( t^n s / n!) {\rm ad}(X_0)^n Y)$ lies in $[L(N)(\l_1), L(N)(\l_1)]$.
This implies that the element $\exp_U ( \sum_{n \geq 1} ( t^n s / n!) {\rm ad}(X_0)^n Y)$ belongs to the group $[N(\l_1),N(\l_1)]$.
Hence appealing to \eqref{u-1} we conclude that, for all $s,t$, 
$$
\phi (  \exp_U ( \sum_{n \geq 1} ( t^n s / n!) {\rm ad}(X_0)^n Y)) = e .
$$
As $U_1$ is an abelian normal subgroup of $U$, from \eqref{u5}, we obtain
\begin{equation}\label{u6}
 \phi( \exp_U (t X_0)) \, \phi( \exp_U (s Y)) \, \phi( \exp_U (t X_0))^{-1}
= \phi(\exp_U ( sY))\phi(\exp_U ( ts [X_0,Y])).  
\end{equation}

Now as $\widetilde{U}_1$ is abelian, using \eqref{u6} and \eqref{u2} we get that
\begin{equation}\label{u7}
  \phi( \exp_U (t X_0)) \, \phi( \exp_U (s Y)) \, \phi( \exp_U (t X_0))^{-1}
= \exp_{\widetilde{U}} ( \sum_l \{\mu_l (s) + \delta_l ( ts) \} \widetilde{Y}_l).
\end{equation}

As $\exp_{\widetilde{U}}$ is a $\l_2$-isomorphism (of $\l_2$-varieties) from $L(\widetilde{U})$ to $\widetilde{U}$ and $\{\, Y_l, \,| 1 \leq l \leq q \}$ is
a $\l_2$-basis for $L(\widetilde{U})(\l_2)$,  
comparing \eqref{u4} and \eqref{u7}, we conclude that, for all $l, \, 1 \leq l \leq q$ and for all $s, t \in \l_1$, 
$$
\delta_l (ts) = \sum_{i,j} a_{ijl} \lambda_i (t) \mu_j(s).
$$

We now apply Lemma \ref{additive-homom} to see $\delta_l =0$, for all $l$. In particular it follows that 
$$\phi( \exp_U (t X_0)) \, \phi( \exp_U ( Y)) \, \phi( \exp_U (t X_0))^{-1}
= \phi(\exp_U ( Y))
$$
for all $t \in \l_1$. In other words $\phi( \exp_U (t X_0))$ commutes with $\phi (v)$, for all $t \in \l_1$.
This completes the proof of Theorem \ref{main-unip}
\hfill$\Box$

\begin{rmk}\label{unip-rmk}{\rm Let $\l_1, \l_2$ be fields such that 
there is an embedding $\sigma : \l_1 \to {\rm M}_n ( \l_2)$ of rings, for some $n$. Let $U$ be a unipotent group over ${\mathbb F_1}$. 
Using $\sigma$ we will now easily construct a unipotent group $\widetilde{U}$ over $\l_2$ and 
an abstract embedding $\eta: U (\l_1) \to \widetilde{U}(\l_2)$. 
We may assume that there is a general linear group, say, ${\rm GL}_m ({\overline\l_1})$ containing  $U$ 
as a unipotent Zariski-closed subgroup which is defined over $\l_1$ (here  ${\rm GL}_m ({\overline\l_1})$ is equipped with its usual $\l_1$-structure
so that the group of $\l_1$-points is ${\rm GL}_m (\l_1)$). 
Now composing the ring embedding 
$\eta_{\sigma} : {\rm M}_m ( \l_1) \to {\rm M}_m ({\rm M}_n (\l_2))$, defined by the assignment $\eta_\sigma ((a_{ij})) : = (\sigma(a_{ij}))$,  $(a_{ij}) \in 
 {\rm M}_m ( \l_1)$, with the obvious $\l_2$-algebra isomorphism ${\rm M}_m ({\rm M}_n (\l_2)) \simeq {\rm M}_{mn} (\l_2)$  we obtain a ring embedding $\eta:{\rm M}_m ( \l_1) \to {\rm M}_{mn} (\l_2)$.  It is immediate that $\eta$ carries the
unipotent elements of ${\rm GL}_m ( \l_1)$ to the unipotent elements in ${\rm GL}_{mn} ( \l_2)$. Thus if $\widetilde{U}:= {\rm Zcl}\,\eta(U (\l_1))$  
then $\widetilde{U}$ is a unipotent group over $\l_2$ and the restriction $\eta : U (\l_1) \to \widetilde{U}(\l_2)$ provides an abstract embedding.
To see an immediate application of the above construction in some specific situation we recall that there is an abstract embedding 
$\q \hookrightarrow \c$ and consequently there is a ring embedding $\q \hookrightarrow {\rm M}_2 (\r)$ (although the only ring homomorphism from $\q$ to $\r$ is the zero homomorphism). Thus, in view of the above observations, it follows that  
if $U$ is a unipotent group over $\q$ then there are unipotent groups $U_1, U_2$ over $\c$ and $\r$, respectively, such that $U (\q)$ embeds abstractly in both the groups $U_1 (\c), U_2 (\r)$.}
\end{rmk}

\begin{rmk}{\rm Let $\l_1, \l_2$ be fields with ${\rm Char} \l_1= {\rm Char}\l_2=0$ and $U, \tilde{U}$ be unipotent groups over $\l_1, \l_2$, respectively.
Let $\phi : U (\l_1) \to \tilde{U}(\l_2)$ be an abstract homomorphism. Considering the conjecture 
due to A. Borel and J. Tits in \cite[\S 8.19]{BT} (see also the conjecture BT in \cite[\S 1]{Ra})  we see that an exact analogue of the conjecture fails to hold if the groups in the domain and the range are replaced by unipotent
groups. However, in view of Theorem \ref{main-unip} and Remark \ref{unip-rmk}  it will be very interesting
to know if a description exists for a general abstract homomorphism $\phi$ between the groups of rational points of unipotent groups, along the same line as the above conjecture,
involving a ring embedding of $\l_1$ in a finite dimensional
algebra over $\l_2$ and an abstract homomorphism from the abelian group $U (\l_1) / [U(\l_1), U(\l_1)]$.}
\end{rmk}

We need the following lemma to draw Corollary \ref{main-unip-cor} from Theorem \ref{main-unip}.

\begin{lemma}\label{divalg} 
Let $\l$ be a field and let $\l^*$ be the multiplicative group of non-zero elements of $\l$. 
Assume that $\l^*$ has a $p$-divisible element which is of infinite order. Let $n$ be an integer. Then  
any ring homomorphism $\rho : \l \to {\rm M}_n (\q)$ is trivial, that is, $\rho(\l) =0$.
Thus,
\begin{enumerate}
\item For an algebraically closed field $\l$,  
any ring homomorphism $\rho : \l \to {\rm M}_n (\q)$ is trivial.
\item Any ring homomorphism $\rho : \r \to {\rm M}_n (\q)$ is trivial.
\item If $l \neq p$ is a prime number then  
 any ring homomorphism $\rho : \Q_l \to {\rm M}_n (\q)$ is trivial.
\end{enumerate}
\end{lemma}  

\noindent{\bf Proof.}
It is enough to show that there is a non-zero element $w$ in the field $\l$ such that $\rho(w) =0$.
Let $t$ be a $p$-divisible element in the multiplicative group $\l^*$ which is of infinite order. 
Clearly $\rho (t)$ is an element in ${\rm GL}_n (\q)$ which is a $p$-divisible element in ${\rm GL}_n (\q)$.    
Appealing to \cite[Theorem 1.5]{Ch1} we see that there is a positive integer $r$ such that $\rho(t)^r$ is a
unipotent matrix in ${\rm GL}_n (\q)$. In particular, as $\rho$ is a ring homomorphism, we have that
$\rho ( (t^r-1)^{2n} ) =0$. As $t$ is not of finite order $t^r-1 \neq 0$. We choose $w = (t^r-1)^{2n}$.

To prove $(1)$ and $(2)$ we need to exhibit an infinite order, $p$-divisible element in $\l^*$, which is immediate when $\l$ is either algebraically closed or $\l= \r$.

As above, to prove (3) we need to exhibit an infinite order, $p$-divisible element in $\Q_l^*$. 
Consider the (multiplicative) group of units in the ring of integers ${\mathbb Z}_l$ of $\Q_l$ and define the subgroup $H_l$ by setting $H_l := 1 + l {\mathbb Z}_l$.
Then observe that $H_l$ is a compact $l$-adic group and $\O (H_l) = l^{\infty}$. Appealing to 
\cite[Theorem 3.3]{Ch1}, as $l \neq p$, we see that the group $H_l$ is $p$-divisible. It is clear that there are at the most
countably many of elements in $H_l$ which are of finite order. On the other hand $H_l$ has uncountably many elements.
This forces $H_l$ to have uncountably many elements which are $p$-divisible and of infinite order.
\hfill$\Box$

\begin{rmk}\label{rmk-divalg}
{\rm We note that the assertion in Lemma \ref{divalg} still holds if $\q$ replaced by a finite extension or if the field $\l$ is replaced by 
a division ring $D$ having the same property as that of the field $\l$, namely, the multiplicative group of invertible elements $D^*$ has a $p$-divisible element of infinite order.}
\end{rmk}

We next prove Theorem \ref{low-card} regarding abstract homomorphism from a unipotent group to an abstract group with smaller cardinality than that of the field over which the unipotent group is defined.

\noindent
{\bf Proof of Theorem \ref{low-card}.}
Before proceeding towards the proof we observe a simple fact on nilpotent
operators on finite dimensional vector spaces over the field $\l$. Let $V$ be a finite dimensional vector space
over $\l$ and let $N \in {\rm End}_{\l}(V)$ be a nilpotent operator of $V$. Let ${\rm GL}(V)$ be the group
of invertible elements in ${\rm End}_{\l}(V)$ and let ${\rm Id}$ denote the identity element.
Then we claim that if $t \in \l$ and $t \neq 0$ then 
\begin{equation}\label{lc1}
(\exp_{\rm{GL}(V)} (N)- {\rm Id})V = (\exp_{\rm{GL}(V)} (tN)- {\rm Id})V.
\end{equation}

To prove \eqref{lc1}, as $t \neq 0$, it is enough to prove that $(\exp_{\rm{GL}(V)} (tN)- {\rm Id})V \subset (\exp_{\rm{GL}(V)} (N)- {\rm Id})V$. Note that, 
for any integer $r >0$ one has $\exp_{\rm{GL}(V)} (rN) = \exp_{\rm{GL}(V)} (N)^r$ and in particular,
\begin{gather}
(\exp_{\rm{GL}(V)} (rN)- {\rm Id})V=\nonumber\\ 
\label{lc2}(\exp_{\rm{GL}(V)} (N)- {\rm Id})( \exp_{\rm{GL}(V)}((r-1)N) + \cdots + {\rm Id})V \subset (\exp_{\rm{GL}(V)} (N)- {\rm Id}) V.
\end{gather}

Let  ${\mathbb G}_a$ be the standard one-dimensional unipotent group over $\l$.
As $N$ is nilpotent, for any $v \in V$ the map $\lambda_v : {\mathbb G}_a (\l) \to V$ defined by $\lambda_v (t) := (\exp_{\rm{GL}(V)} (tN)- {\rm Id}) v$, $t \in 
\l$, is algebraic (and defined over $\l$). 
Further, by the inclusion as in \eqref{lc2}, it follows that $\lambda_v (r) \in (\exp_{\rm{GL}(V)} (N)- {\rm Id})V$ for all $r \in {\mathbb N}$. Hence, as ${\mathbb N}$ is Zariski-dense in $\overline{\l}$, we conclude that 
$\lambda_v (t) \in (\exp_{\rm{GL}(V)} (N)- {\rm Id})V$ for all $t \in \l$. This proves the claim.

We will now get back to the proof of theorem. 
We will prove the theorem by induction on the dimension of $U$.
In the case when ${\rm dim}\,U =1$ the group $U (\l)$ is abelian and the statement follows trivially.
We next assume that the statement holds for all unipotent groups $U'$ defined over $\l$ such that  ${\rm dim}\,U' < {\rm dim}\,U$.
We will prove that the statement holds for $U$. As $U$ is $\l$-split there is a normal $\l$-subgroup $U_1$ of $U$ such that
$U/U_1$ is $\l$-isomorphic to the one-dimensional additive group ${\mathbb G}_a$ (defined over $\l$).
This says that there is a element $X_0 \in L(U)(\l)$ such that if $E$ denotes the Zariski closure of the one-parameter subgroup
$\{\exp_U (tX_0) \, | \, t \in \l \}$ then 
$U(\l) = E(\l) U_1 (\l)$. Now as ${\rm dim}\,U_1 < {\rm dim}\,U$, by induction hypothesis, the image $\phi (U_1 (\l))$
is an abelian subgroup of $H$. Hence it is enough to prove that if $g \in E(\mathbb{F})$ and $g$ is not the identity element of $U(\l)$
then the element $\phi (g)$ commutes with $\phi (U_1 (\l))$. 
As the cardinality of $\l$ is strictly greater than that of $H$ there is a $t_0 \neq 0$ such
that $\phi (\exp_U (t_0X_0)) = e$, where $e$ denotes the identity element of $H$. Setting $X_1 : = t_0X_0$ we see that 
$E(\mathbb{F}) = \{\exp_U (tX_1) \, | \, t \in \l \}$. Let $ s \in \l$ such that $g = \exp_U (s X_1)$. As $g$ is not
the identity element of $U(\l)$ it is immediate that $s \neq 0$.
We apply \eqref{lc1} to see that
\begin{equation}\label{lc3}
({\rm Ad} (\exp_U(X_1)) - {\rm Id})L(U_1)(\l) =  ({\rm Ad} (\exp_U(sX_1)) - {\rm Id})L(U_1)(\l).
\end{equation}

It follows from \cite[Proposition 1.2]{F} that $[U_1,U_1] (\l) = [U_1 (\l),U_1 (\l)]$. In view of this fact,
exponentiating both sides of the equation \eqref{lc3} we get that, for any $Y \in L(U_1) (\l)$ there is a $h \in [U_1 (\l),U_1 (\l)]$
and $Z \in  L(U_1) (\l)$ such that
\begin{equation}\label{lc4}
[g, \exp_U(Y)] = [\exp_U (X_1), \exp_U(Z)] h.
\end{equation}

As $\phi (U_1 (\l))$ is abelian $\phi (h) =e$ and as $\phi (\exp_U(X_1)) = e$ it follows from the equation \eqref{lc4} that
$\phi(g)$ commutes with $\phi (\exp_U(Y))$ for all $Y \in L(U_1)(\l)$. In particular, $\phi(g)$ commutes with all of $\phi (U_1 (\l))$.
This completes the proof of Theorem \ref{low-card}.
\hfill$\Box$

We next turn our attention to images of solvable algebraic groups under abstract homomorphisms. 
We first observe a short but convenient lemma.

\begin{lemma}\label{unip-trivial} 
Let $\l$ be a field with ${\rm Char}\l =0$ and $G$ be an algebraic group over $\l$. Let $g, u \in G (\l)$ be 
unipotent elements and $m \neq 1$ be an integer such that $g u g^{-1} = u^m$. Then $u =e$.
\end{lemma}

\noindent{\bf Proof.}  
Recall that  $\exp_G : {\mathcal N}_G \to {\mathcal U}_G $ is
a $\l$-isomorphism between the variety of nilpotent elements ${\mathcal N}_G$ of $L(G)$ and the variety of unipotent elements ${\mathcal U}_G$ of $G$. Thus
there is a nilpotent element $X \in L(G)( \l)$ such that $\exp (X) =u$. Now the hypothesis
in the lemma implies that  
$$\exp_G ({\rm Ad}(g) X) = \exp_G (mX).
$$
In particular, as $\exp_G : {\mathcal N}_G \to {\mathcal U}_G $ is a isomorphism of $\l$-varieties, it follows that
${\rm Ad}(g) X = mX$. As ${\rm Ad}(g)$ is a unipotent operator of $L(G)$ and as $m \neq 1$ we conclude that
$X=0$. In other words, we have, $u=e$.
\hfill$\Box$

\begin{prop}\label{prop-solv} 
Let $G$ be a linear algebraic group defined over $\k$. Let $\l $ be a field, $B$ be a (Zariski-)connected 
solvable group over $\l$ and $\phi : B(\l) \to G(\k)$ be an abstract homomorphism. 
Let $T$ be a maximal $\l$-torus of $B$.
Let $M \subset R_u (B)$ be a normal $\l$-subgroup of $B$ such that
$B= Z_B (T) M$ and $M \cap Z_B(T) =e$. Then $\phi (M(\l)) =e$ in the following three cases.

\begin{enumerate}
\item $\l$ is algebraically closed field with ${\rm Char}\l =0$.
\item $\l= \r$. 
\item $l$ is a prime with $l \neq p$, $\l= \Q_l$ and $B$ is $\Q_l$-split.
\end{enumerate}
\end{prop}

\noindent{\bf Proof.} 
Let $B' =  B / [R_u (B),R_u (B)]$. It is clear that $R_u (B')$ is abelian. 
Let $\pi : B \to B'$ be the natural quotient homomorphism. Let $M \subset R_u (B)$ be as in the statement of the proposition.
We set $T' := \pi (T)$ and $M' := \pi (M)$. Clearly $T'$ is a maximal $\l$-torus of $B'$ and moreover, $B' = Z_{B'} (T') M'$.  
It is not difficult to see that $\pi (  Z_B (T)) = Z_{B'} (T')$ and $Z_{B'} (T') \cap M' =e$.

As $R_u(B)$ is unipotent, appealing to Corollary \ref{main-unip-cor}, we see that, when $\l$ is
either an algebraically closed field or the field of 
real numbers $\r$ or the $l$-adic field $\Q_l$, $l \neq p$, the image
$\phi (R_u(B) (\l))$ is an abelian subgroup
of $G(\k)$. Thus the homomorphism $\phi$ factors through and abstract homomorphism 
$\phi' : B(\l)/ [R_u(B),R_u(B)] (\l) \to G(\k)$.
It follows from  \cite[Lemma 2.3]{F} that the groups $B'(\l)$ and $B(\l)/ [R_u(B),R_u(B)] (\l)$ are isomorphic. Thus to prove that 
$\phi (M(\l)) =e$ it is enough to show that $\phi': (M' (\l)) =e$. 
Consequently, without loss of generality, we may further assume that $R_u (B)$ is abelian. 

We now consider three cases according as $\l$ is either algebraically closed or $\l = \r$ or 
$\l = \Q_l, \, l \neq p$ and $B$ is $\Q_l$-split. 

{\it Case 1:} We assume that $\l$ is algebraically closed.  
We decompose $L(M)(\l)$ 
as $L(M)(\l) =  L(M)_{\chi_1}(\l) + \cdots + L(M)_{\chi_k}(\l)$ 
where $\chi_i, \, 1 \leq i \leq k$ are distinct (multiplicative) algebraic characters on $T(\l)$
and $L(M)_{\chi_i}(\l), \, 1 \leq i \leq k$ is the non-zero subspace of $L(M)(\l)$ on which an element $z \in T(\l)$ 
acts by multiplication by $\chi_i (z)$. Since 
$T (\l)$ does not fix non-zero elements in $L(M) (\l)$
it follows that $\chi_i \neq 1$ for all $1 \leq i \leq k$. To prove that $\phi (M(\l))=e$ it is enough to show that
$\phi (\exp_{R_u(B)} ( L(M)_{\chi_i}(\l)))=e$, for all $1 \leq i \leq k$. We consider  $ X \in L(M)_{\chi_i}(\l)$. As $\chi_i \neq 1$ there is a
$t \in T(\l)$ such that  $\chi_i (t) = -1$. Thus $t \exp_{R_u(B)} (X) t^{-1} = \exp_{R_u(B)} (X)^{-1}$. Applying $\phi$ we obtain
$\phi(t) \phi(\exp_{R_u(B)} (X)) \phi(t)^{-1} = \phi(\exp_{R_u(B)} (X))^{-1}$. 
As both $t,\exp_{R_u(B)} (X)$ are divisible elements it follows that both $\phi(t), \phi(\exp_{R_u(B)} (X))$ are unipotent elements in $G(\k)$.
Thus, appealing to Lemma \ref{unip-trivial} we get that $\phi(\exp_{R_u(B)} (X))=e$. 

{\it Case 2: } We assume that $\l = \r$. We first write $T = T_a T_s$
where $T_a$ is the $\r$-anisotropic part of $T$ and $T_s$ is the $\r$-split part of $T$.
As $T_s$ is $\r$-split we 
we decompose $L(M)(\r)$ 
as $L(M) (\r) =  L(M)_{\chi_1}(\r) + \cdots + L(M)_{\chi_k}(\r)$ where $\chi_i, \, 1 \leq i \leq k$ are distinct (multiplicative) algebraic characters on $T_s$
and $L(M)_{\chi_i}(\r), \, 1 \leq i \leq k$ is the subspace of $L(M)(\r)$ on which an element $z \in T_s (\r)$ acts by multiplication by $\chi_i(z)$.
As in Case 1, it is enough to prove that 
$\phi (\exp_{R_u(B)} ( L(M)_{\chi_i}(\r)))=e$, for all $1 \leq i \leq k$. We need to consider two cases. 

First assume that $\chi_i \neq 1$. Recall that $T_s(\r)^*$ denotes the connected component of $T_s(\r)$ in the topology induced by $\r$. As 
$\chi_i \neq 1$ we may choose an integer $m \neq 1$ and $t \in T_s (\r)^*$ such that $\chi_i (t) = m$. 
Thus ${\rm Ad}(t) v = mv$, for all $v \in L(M)_{\chi_i}(\r)$. Hence $t \exp_{R_u(B)} (v) t^{-1} = (\exp_{R_u(B)} (v))^m$. 
Recall that  $T_s (\r)^*$ is a divisible group.
Thus both $t,\exp_{R_u(B)} (X)$ are divisible elements and it follows that both $\phi(t), \phi(\exp_{R_u(B)} (X))$ are unipotent elements in $G(\k)$.
Appealing to Lemma \ref{unip-trivial} we get that $\phi(\exp_{R_u(B)} (X))=e$.

We next assume that $\chi_i =1$. Thus $T_s (\r)$ acts trivially on $L(M)_{\chi_i}(\r)$. 
It is clear that $L(M)_{\chi_i}(\r)$ remains invariant under the action of the compact part $T_a (\r)$.
Since no non-zero vector in $L(M)_{\chi_i}(\r)$  remains fixed by all the elements in 
$T(\r)$ it is immediate that $T_a( \r)$ has no non-zero fixed point in $L(M)_{\chi_i}(\r)$. Since $T_a (\r)$ is compact abelian we may further
decompose $L(M)_{\chi_i}(\r)$ into  $T_a( \r)$-irreducible components as 
$L(M)_{\chi_i}(\r) = V_1 + \cdots + V_l$. To prove $\phi (\exp_{R_u(B)} ( L(M)_{\chi_i}(\r)))=e$ it is now enough to show that  
$\phi (\exp (V_j)) = e$, for each $ 1 \leq j \leq l$. 
As $T_a (\r)$ is compact, connected and abelian it is a direct product of circle groups and hence one of the direct factor, say $C$, of 
$T_a (\r)$ acts non-trivially on $V_j$. 
In particular there is an element $c \in C$
such that ${\rm Ad} (c) v = -v$, for all $v \in V_j$. Hence $c \exp_{R_u(B)} (v) c^{-1} = \exp_{R_u(B)} (v)^{-1}$, for all $v \in V_i$
Now as $T_a(\r)$ is a divisible group it follows that $\phi (c)$ is divisible element in $G(\k)$.
We use (3) of Theorem \ref{key1} to see that
 $\phi(c)$ is a unipotent element in $G(\k)$. Analogously,
it is clear that $\phi (\exp_{R_u(B)} (v))$ is unipotent, for all $v \in V_j$. 
We now use Lemma \ref{unip-trivial} to see that $\phi (\exp_{R_u(B)} (v)) = e$, for all $v \in V_j$

{\it Case 3: } We assume that $\l = \Q_l$ and $B$ is $\Q_l$-split.
As $B$ is $\Q_l$-split the torus $T$ is $\Q_l$-split. For an integer $r$ we define
$W_r := 1 + l^{2+r} {\mathbb Z}_l \subset \Q_l^*$. Let ${\mathbb G}_m$ be the standard one-dimensional $\l$-split torus and
and let ${\mathbb G}^n_m$ be the product
of $n$ copies of  ${\mathbb G}_m$, where $n := {\dim T}$ .
We first fix a $\Q_l$-isomorphism $\lambda :{\rm {\mathbb G}}^n_m \to T$ and define
$T(W_r):= \lambda ( W_r^n)$.  It is easy to see that $T(W_r)$ is a pro-$l$ group and in particular $\O (T(W_r)) = l^\infty$. 
As $l \neq p$ using \cite[Lemma 3.2]{Ch1} the group $T(W_r)$ is $p$-divisible.
We now appeal to Proposition \ref{p-divisible} to conclude that ${\rm Zcl}\, \phi(T(W_1)) = F U'$ where $F$ is a finite abelian group,
$U'$ is a unipotent group and $F$ normalises $U'$. Further we may assume that $F$ is an abstract quotient of the group $T(W_1)$.
Hence using Theorem \ref{profinite} we see that $\O(F) = l^s$ for some integer $s$. Thus it follows that
$\phi(P_{l^s} (T(W_1))) \subset U' (\k)$, where $P_{l^s}: T(W_1) \to T(W_1)$ is the $l^s$-th power map. In particular it follows that  if $ z \in  T(W_s)$
then $\phi (z)$ is a unipotent element in $G (\k)$. 

We next give an argument similar to that of the previous two cases. 
As $T$ is $\Q_l$-split, we decompose $L(M)(\Q_l)$ 
as $L(M)(\Q_l) =  L(M)_{\chi_1}(\Q_l) + \cdots + L(M)_{\chi_k}(\Q_l)$ where $\chi_i, \, 1 \leq i \leq k$ 
are distinct (multiplicative) algebraic characters on $T$
and $L(M)_{\chi_i}(\Q_l), \, 1 \leq i \leq k$ is the subspace of 
$L(M)(\Q_l)$ on which an element $z \in T(\Q_l)$ acts by multiplication by $\chi_i (z)$. Since 
$T (\Q_l)$ does not fix non-zero elements in $L(M) (\Q_l)$
it follows that $\chi_i \neq 1$ for all $1 \leq i \leq k$.
As before it is enough to prove that $\phi (\exp_{R_u (B)} (L(M)_{\chi_i}(\Q_l)) = e$.  
Since $\chi_i \neq 1$ it is not difficult to see that there is an integer $m \neq 1$ and $t \in T(W_s)$ such that $\chi(t) =m$.
Thus ${\rm Ad}(t) v = mv$, for all $v \in L(M)_{\chi_i}(\r)$. Hence $t \exp_{R_u(B)} (v) t^{-1} = (\exp_{R_u(B)} (v))^m$. As both $\phi (t)$ and $\phi (\exp_{R_u(B)} (v))$
are unipotent elements in $G(\k)$ we appeal to Lemma \ref{unip-trivial} to conclude that 
$\phi (\exp_{R_u(B)} (v)) = e$, for all $v \in L(M)_{\chi_i}(\Q_l)$.
This completes the proof.
\hfill$\Box$

We need to recall a result is due to A. Borel and T. Springer which ensures existence of a subgroup $M$ of $B$ as in Proposition \ref{prop-solv}.

\begin{thm}\label{bo-sp}{\rm (Borel-Springer \cite[Lemma 9.7, p. 487]{BS})}
Let $B$ be a connected solvable group over $\l$ and $T$ be a $\l$-torus of $G$. Let $R_u (B)$ be the unipotent radical
of $B$. Then $B$ has a normal $\l$-subgroup $M$ contained in $R_u (B)$ such that $B= Z_B (T) M$ and $M \cap Z_B(T) =e$.
\end{thm}

We now apply Proposition \ref{prop-solv} and Theorem \ref{bo-sp} to deduce Theorem \ref{thm-solv}.

\noindent
{\bf Proof of Theorem \ref{thm-solv}.}
Using  Proposition \ref{prop-solv} we first prove that $\phi (B(\l))$ is abelian in all the three cases.
Let $T$ be a maximal $\l$-torus of $B$.
Let $M$ be the normal $\l$-closed subgroup of $B$ as given in Theorem \ref{bo-sp}. 
Then $B = Z_B(T)M$ and moreover $B(\l) = Z_B(T) (\l) M(\l)$. 
Appealing to Proposition \ref{prop-solv} we see that $\phi(B(\l)) = \phi (Z_B(T) (\l))$. 
Hence it is enough to prove that $\phi (Z_B(T) (\l))$ is abelian.
One may further write 
$Z_B (T) = T Z_{R_u (B)} (T)$ and in particular, 
$Z_B(T) (\l) = T(\l)Z_{R_u (B)} (T)(\l)$. 
Now we apply Corollary \ref{main-unip-cor} to see that $\phi (Z_{R_u (B)} (T)(\l))$ is abelian (and unipotent).
As both $\phi (Z_{R_u (B)} (T)(\l))$, $\phi (T(\l))$ are abelian 
subgroups and as 
$Z_{R_u (B)} (T)(\l)$ commutes with $T(\l)$ 
it is clear that $\phi (Z_B(T) (\l))$ is abelian. Thus in all the three cases 
$\phi (B(\l))$ is abelian.

It now remains to prove (1) and (2).
To prove (1) observe that, as $T(\l)$ is divisible,
the group $\phi (T(\l))$ is unipotent and in particular $\phi (B(\l))$ is unipotent.
We prove (2) analogously, by noting that $T(\r)^*$ is a divisible group.

This completes the proof of Theorem \ref{thm-solv}.
\hfill$\Box$

\section{Homomorphisms from non-solvable groups}
This section is devoted to studying abstract group homomorphisms from the groups of rational points of 
non-solvable algebraic groups. 
One of the goals in this section is to give proofs of  Theorem \ref{real-image}, Theorem \ref{real},
Theorem \ref{alg-closed}, Theorem \ref{isotropic}, Theorem \ref{sl(D)} and Theorem \ref{final-l-adic}. 
As mentioned in \S 1, we fix a prime number $p$ and $\k$ will always denote a number field or a finite a extension of $\q$.
We begin this section by formulating  general results in Proposition \ref{borel-tits-app} and in Corollary \ref{borel-tits-app-cor-1} concerning abstract homomorphisms 
between groups of rational points over algebraic groups over general fields, which are devised with the aim of proving Theorem \ref{real-image}.  As immediate consequences we also obtain
Corollary \ref{borel-tits-app-cor-2}, Corollary \ref{borel-tits-app-cor-3} and  Corollary \ref{borel-tits-app-cor-4}.
The proof of Proposition \ref{borel-tits-app} follows applying the landmark Theorem \ref{borel-tits} due to Borel and Tits, and the celebrated result (4) of Theorem \ref{kneser-tits} due to Tits.

\begin{prop}\label{borel-tits-app} 
Let $\l_1, \l_2$ be  fields such that $\l_1$ is infinite and $\l_2$ is perfect.
Let $H$ and $G$ be algebraic groups over $\l_1$ and $\l_2$, respectively. Moreover assume that
$H$ is $\l_1$-isotropic and absolutely simple. Let $\phi : H (\l_1)^+ \to G (\l_2)$ be 
an abstract group homomorphism. If $\phi$ is non-trivial then there exist 
\begin{enumerate}
\item a finite extension
$\k_2$ of $\l_2$ and a field homomorphism $\alpha : \l_1 \to  \k_2$ and

\item a surjective algebraic group homomorphism 
$\mu : {\rm Zcl}\,\phi (H (\l_1)^+) \to {\mathcal R}_{\k_2/\l_2}( \! {^\alpha \! H} / Z(\! {^\alpha \! H} ))$ over $\l_2$
\end{enumerate}
such that $\pi   \circ \mu \circ \phi = \delta \circ \alpha^0$,
where $\delta : \! {^\alpha \! H}  \to \! {^\alpha \! H} / Z(\! {^\alpha \! H} )$ is the 
homomorphism from the group to the quotient, 
and $ \pi : {\mathcal R}_{\k_2/\l_2}( \! {^\alpha \! H} / Z(\! {^\alpha \! H} ))
\to  \! {^\alpha \! H} / Z(\! {^\alpha \! H} )$ is the natural projection over $\k_2$ which induces
the isomorphism between ${\mathcal R}_{\k_2/\l_2}( \! {^\alpha \! H} / Z(\! {^\alpha \! H} )) (\l_2)$ and
$ \! {^\alpha \! H} / Z(\! {^\alpha \! H} ) (\k_2)$ when $\pi$ is restricted to the former.
\end{prop} 

\noindent
{\bf Proof.} 
Let $A := {\rm Zcl}\, (\phi ( H^+ (\l_1)))$. In the case when $ {\rm Char \l_1}=0$ we see in the following elementary way that  $A$ Zariski-connected.
Let $A^0$ be the Zariski connected component of the identity element.
As $H^+ (\l_1)$ is generated by unipotent elements of $H(\l_1)$ it is enough to show that,
for any unipotent element $u \in H(\l_1)$ the element $\phi (u)$ lies in $A^0$.
Note that, as $u$ is unipotent and $ {\rm Char \l_1}=0$, for any integer $n$ there is a $v \in H(\l_1)$ such that $u =v^n$.
Thus $\phi (u)A^0$ has a $n$-th root for any integer $n$ in the finite group $A/A^0$.
This implies that $\phi (u) \in A^0$. Hence $A = A^0$. 
Now the general case, when $\l_1$ is assumed only to be infinite, follows from \cite[Proposition 7.2 (i), p. 539]{BT}.

We now assume that $\phi$ is non-trivial and then arrive at the conclusion. 
Clearly the group $A$ is defined over $\l_2$ and, moreover, as $\l_2$ is perfect it follows that $R_u (A)$ is defined over $\l_2$.
We now consider the reductive algebraic group $A/ R_u (A)$ which has a $\l$-structure induced from that of $A$. 
Using (4) Theorem \ref{kneser-tits}, we have $[H(\l_1)^+ , H(\l_1)^+ ] = H(\l_1)^+$.
Thus non-triviality of $\phi$ implies that $A$ is non-solvable. In particular ${\rm dim}\, A/ R_u (A)  > 0$. It is now easy to see that $[A/ R_u (A), A/ R_u (A)] = A/ R_u (A)$.
This follows from the fact that $[H(\l_1)^+ , H(\l_1)^+ ] = H(\l_1)^+$ and that $ \nu \circ \phi ( H(\l_1)^+)$ is Zariski dense in $A/ R_u (A)$, where $\nu : A \to A/ R_u (A)$ is the quotient homomorphism over $\l_2$. 
Thus $A/ R_u (A)$ is in fact a 
semisimple group over $\l_2$. Let $S$ be a non-trivial $\l_2$-simple factor of $A/ R_u (A)$.
It is a standard fact that, taking quotient of $A/ R_u (A)$ by the product of all the
$\l_2$-simple factors of $A/ R_u (A)$ (if any) other than $S$, one would obtain a surjective algebraic group
homomorphism over $\l_2$, say, $\widetilde{\lambda} : A/ R_u (A) \to S/ Z(S)$. We now consider the surjective algebraic group
homomorphism over $\l_2$, say, $\lambda : A  \to S/ Z(S)$ defined by $\lambda : = \widetilde{\lambda} \circ \nu$.
As $S/ Z(S)$ is $\l_2$-simple there exists
a finite extension $\k_2$ of $\l_2$ and an absolutely simple and adjoint $\k_2$-algebraic group $B$
such that $S / Z(S) = {\mathcal R}_{\k_2/ \l_2} (B)$. Let $\pi' : {\mathcal R}_{\k_2/ \l_2} (B)
\to B$ be the projection homomorphism which is defined over $\k_2$ such that $\pi' : {\mathcal R}_{\k_2/ \l_2} (B) (\l_2) \to B (\k_2)$ is an isomorphism of abstract groups. 
We consider the abstract homomorphism $\pi' \circ \lambda \circ \phi : H(\l_1)^+ \to B (\k_2)$, and
observe that, as
$\lambda \circ \phi ( H (\l_1)^+)$ is Zariski-dense in $S/Z(S)$ so is 
$\pi' \circ \lambda \circ \phi ( H (\l_1)^+)$ in $B$.
We now appeal to the Theorem \ref{borel-tits} to get a field homomorphism $\alpha : \l_1 \to \k_2$ 
and a $\k_2$-isogeny $\gamma : \! {^\alpha \! H} \to B$ such that $\pi' \circ \lambda \circ \phi = \gamma \circ \alpha^0$.
Recall that $\delta : \! {^\alpha \! H}  \to \! {^\alpha \! H} / Z(\! {^\alpha \! H} )$ denotes the usual quotient homomorphism.
As $B$ is adjoint there is a $\k_2$-isomorphism, say,
$\bar{\gamma} : \! {^\alpha \! H} / Z(\! {^\alpha \! H} )\to B$
such that $ \bar{\gamma} \circ \delta = \gamma$.
Let ${\mathcal R}_{\k_2/\l_2} (\bar{\gamma}) :
{\mathcal R}_{\k_2/\l_2} ( \! {^\alpha \! H} / Z(\! {^\alpha \! H} )) \to {\mathcal R}_{\k_2/\l_2} (B)$
be the $\l_2$-isomorphism such that $\pi \circ \bar{\gamma} = {\mathcal R}_{\k_2/\l_2} (\bar{\gamma}) \circ \pi'$.
We now compete the proof by setting $\mu := {\mathcal R}_{\k_2/\l_2} (\bar{\gamma})^{-1} \circ \lambda \circ \phi$.
\hfill$\Box$

The next result gives a necessary restriction on the dimensions on groups involved for the existence of non-trivial
abstract homomorphisms.   

\begin{cor}\label{borel-tits-app-cor-1} 
Let $\l_1, \l_2$ be  fields such that $\l_1$ is infinite and $\l_2$ is perfect.
Let $H$ and $G$ be algebraic groups over $\l_1$ and $\l_2$, respectively. Assume that
$H$ is $\l_1$-isotropic and absolutely simple. Further assume that there exists a non-trivial  abstract group homomorphism $\phi : H (\l_1)^+ \to G (\l_2)$.
\begin{enumerate}
\item Then ${\rm dim}\, G/R(G) \geq {\rm dim}\, H$.
\item Moreover, if $\l_1$ does not embed in $\l_2$ then ${\rm dim}\, G/R(G) \geq 2{\rm dim}\, H$.
\end{enumerate}
\end{cor}

\noindent
{\bf Proof.} We use Proposition \ref{borel-tits-app} and the notation involved therein.
Using Proposition \ref{borel-tits-app} we see that there exist a finite extension
$\k_2$ of $\l_2$ and a field homomorphism $\alpha : \l_1 \to  \k_2$ and a surjective algebraic group homomorphism 
$\mu : {\rm Zcl}\,\phi (H (\l_1)^+) \to {\mathcal R}_{\k_2/\l_2}( \! {^\alpha \! H} / Z(\! {^\alpha \! H} ))$ over $\l_2$
such that $\pi   \circ \mu \circ \phi = \delta \circ \alpha^0$,
where $\delta : \! {^\alpha \! H}  \to \! {^\alpha \! H} / Z(\! {^\alpha \! H} )$ is the 
homomorphism from the group to the quotient.
As $\mu$ is surjective the existence of a non-trivial abstract group homomorphism $\phi : H (\l_1)^+ \to G (\l_2)$ implies 
\begin{equation}\label{borel-tits-app-cor-1-eq-1}
{\rm dim}\, G \geq {\rm \dim} {\rm Zcl}\,\phi (H (\l_1)^+ ) \geq {\rm dim}\, {\mathcal R}_{\k_2/\l_2}( \! {^\alpha \! H} / Z(\! {^\alpha \! H} )) = [\k_2 : \l_2] {\rm dim } H \geq {\rm dim }\, H.
\end{equation}

Now from (4) Theorem \ref{kneser-tits}, we have $[H(\l_1)^+ , H(\l_1)^+ ] = H(\l_1)^+$. This in turn implies that 
the abstract group homomorphism $\phi : H (\l_1)^+ \to G (\l_2)$ is non-trivial if and only if the abstract group homomorphism $\theta \circ \phi : H (\l_1)^+ \to G/ R(G) (\l_2)$ is non-trivial, where
$R(G)$ is the solvable radical of $G$ and $\theta: G \to G/R(G)$ is the natural quotient homomorphism. Hence replacing $G$ by $G/R(G)$ in \eqref{borel-tits-app-cor-1-eq-1}
we see that the existence of a non-trivial abstract group homomorphism $\phi : H (\l_1)^+ \to G (\l_2)$ implies 
$$
{\rm dim}\, G/ R(G) \geq  [\k_2 : \l_2] {\rm dim } H \geq {\rm dim }\, H.
$$
Thus (1) is proved. The proof of (2) now follows immediately from the above inequalities by noting that if $\l_1$ does not embed in $\l_2$ then $[\k_2 : \l_2] \geq 2$.
\hfill$\Box$

In the following corollary we show that equality of dimensions of the groups in the  domain and in the range gives strong conclusions on the structure of non-trivial abstract homomorphism which, in particular, implies the Zariski-density of the image.

\begin{cor}\label{borel-tits-app-cor-2} 
Let $\l_1, \l_2$ be  fields such that $\l_1$ is infinite and $\l_2$ is perfect.
Let $H$ and $G$ be connected algebraic groups over $\l_1$ and $\l_2$, respectively. Assume that
$H$ is $\l_1$-isotropic, absolutely simple and that ${\rm dim}\, G = {\rm dim}\, H$. If 
there exists a non-trivial  abstract group homomorphism $\phi : H (\l_1)^+ \to G (\l_2)$ then 
\begin{enumerate}
\item ${\rm Zcl}\, (\phi ( H^+ (\l_1))) = G$,
\item $G$ is absolutely simple.
\end{enumerate}
Consequently, if $H$ is simply connected or $Z(G) =e$,  by Theorem \ref{borel-tits}, there exists a unique field homomorphism $\alpha : \l_1 \to \l_2$,
an isogeny $\beta :  {^\alpha \! H} \to G $ of $\l_2$-algebraic groups such that 
$$
\phi (g) = \beta (\alpha^0 (g)), \,\, for \,\, all \,\, g \in G_1( \l_1).
$$
\end{cor}

\noindent
{\bf Proof.} As above we use Proposition \ref{borel-tits-app} and the notation involved therein.
As before, appealing to Proposition \ref{borel-tits-app} we see that there exist a finite extension
$\k_2$ of $\l_2$ and a field homomorphism $\alpha : \l_1 \to  \k_2$ and a surjective algebraic group homomorphism 
$\mu : {\rm Zcl}\,\phi (H (\l_1)^+) \to {\mathcal R}_{\k_2/\l_2}( \! {^\alpha \! H} / Z(\! {^\alpha \! H} ))$ over $\l_2$
such that $\pi   \circ \mu \circ \phi = \delta \circ \alpha^0$,
where $\delta : \! {^\alpha \! H}  \to \! {^\alpha \! H} / Z(\! {^\alpha \! H} )$ is the 
homomorphism from the group to the quotient.
Now in view of \eqref{borel-tits-app-cor-1-eq-1} and the hypothesis
${\rm dim}\, G = {\rm dim}\, H$ it follows that ${\rm Zcl}\, (\phi ( H^+ (\l_1))) = G$ and $[\k_2: \l_2] =1$.
In particular, ${\mathcal R}_{\k_2/\l_2}( \! {^\alpha \! H} / Z(\! {^\alpha \! H} ))= \! {^\alpha \! H} / Z(\! {^\alpha \! H} )$ and  
$\mu : G \to  \! {^\alpha \! H} / Z(\! {^\alpha \! H} )$ is a surjective homomorphism 
of algebraic groups over $\l_2$. As ${\rm dim}\, \! {^\alpha \! H} / Z(\! {^\alpha \! H} ) = {\rm dim}\, H = {\rm dim}\, G$ it follows that ${\rm dim} \, {\rm ker}\mu =0$. In particular,
${\rm ker}\mu$ is a finite central subgroup of $G$. Thus, as $\! {^\alpha \! H} / Z(\! {^\alpha \! H} )$ is absolutely simple, so is $G$. 
This completes the proof of (1) and (2).
\hfill$\Box$

\begin{rmk}\label{rmk-Che}
{\rm In \cite[Theorem 1]{Che}, when $|\l_1| >5$ and $\l_2$ arbitrary,
all non-trivial abstract homomorphism $\phi : {\rm SL}_2 (\l_1) \to  {\rm GL}_2 (\l_2)$ are classified, and
it is proved that, given any such $\phi$, there exists
a field homomorphism $\alpha : \l_1 \to \l_2$ and $A \in {\rm SL}_2 (\l_2)$ such that $\phi (X) = A {^\alpha \! X} A^{-1}$,
for all $X \in {\rm SL}_2 (\l_1)$,
where ${^\alpha \! X}$ denotes the matrix obtained by applying $\alpha$ to all the entries in $X$.
In \cite[Corollary 2.6]{Che} it is also proved that, if $\l_1, \l_2$ are arbitrary fields with 
$\l_1$ infinite then the image of any non-trivial homomorphism
from ${\rm SL}_2 (\l_1)$ to ${\rm SL}_2 (\l_2)$ is Zariski dense.
By the well-known fact (see \cite[Lemma 2, p. 377]{J}),  as $|\l_1| >5$, the commutator 
$[{\rm SL}_2 (\l_1),{\rm SL}_2 (\l_1)]$ coincides with ${\rm SL}_2 (\l_1)$, and in particular,  $\phi ( {\rm SL}_2 (\l_1)) \subset {\rm SL}_2 (\l_2)$.
It is also well-known that ${\rm SL}_2 (\l_1) = {\rm SL}_2 (\l_1)^+$; see \cite[Lemma 1, p. 376]{J}. In view of the above two facts the Corollary \ref{borel-tits-app-cor-2} generalises  \cite[Theorem 1, Corollary 2.6]{Che} in the case when $\l_1$ is infinite and $\l_2$ is perfect.} 
\end{rmk}

We next record an immediate consequence of the above result on the non-trivial abstract endomorphisms.

\begin{cor}\label{borel-tits-app-cor-3} 
Let $\l_1$ be  a field such that $\l_1$ is infinite and  perfect and let $H$ be an algebraic group over $\l_1$ such that 
$H$ is $\l_1$-isotropic, absolutely simple. Moreover, assume that $H$ is either simply connected or $Z(H) =e$. Let 
$\phi : H (\l_1)^+ \to H(\l_1)^+$ be a non-trivial abstract group homomorphism. Then 
there exists a unique field homomorphism $\alpha : \l_1 \to \l_1$,
an isomorphism $\beta :  {^\alpha \! H} \to H $ of $\l_1$-algebraic groups such that 
$$
\phi (g) = \beta (\alpha^0 (g)), \,\, for \,\, all \,\, g \in H( \l_1)^+ .
$$
Consequently, if $\l_1$ is a finite extension of $\Q$ then $\phi$ is an isomorphism, and moreover, if $\l_1 = \Q$ then
$\phi$ is the restriction of an algebraic isomorphism $\psi : H \to H$ which is defined over $\Q$. 
\end{cor}

\noindent
{\bf Proof.} This follows immediately from Corollary \ref{borel-tits-app-cor-2} and the fact
that if $H$ is simply connected or adjoint then so is ${^\alpha \! H}$ and consequently, an
isogeny $\beta : {^\alpha \! H} \to H$ is an isomorphism.
\hfill$\Box$    

\begin{cor}\label{borel-tits-app-cor-4} 
Let $\l_1, \l_2$ be  fields such that $\l_1$ is infinite and $\l_2$ is perfect.
Assume further that $\l_1$ can not be embedded in any finite extension of $\l_2$. Let $H$ be a $\l_1$-simple $\l_1$-isotropic group and $G$ be any algebraic group over $\l_2$. Then any abstract homomorphism $\phi : H (\l_1)^+ \to G ( \l_2)$ is trivial
\end{cor}

\noindent
{\bf Proof.} The proof follows immediately from Proposition \ref{borel-tits-app} and the fact that 
there is a finite separable extension $\k_1$ of $\l_1$ and an absolutely simple
group $E$ over $\k_1$ such that $H (\l_1)^+ \simeq E (\k_1)^+$.
\hfill$\Box$

\noindent{\bf Proof of Theorem \ref{real-image}.}
To prove the first statement in (1) we apply Proposition \ref{borel-tits-app}
by setting $\l_1 := \l$ and $\l_2 := \c$. The second statement in (1) now follows from (1) of Corollary \ref{borel-tits-app-cor-1}.

To prove the first statement of (2) we again apply Proposition \ref{borel-tits-app}, as above.
We set $\l_1 := \l$ and $\l_2 := \r$. 
Then, as $\l$ has no embedding in $\r$, the proof is completed by the immediate conclusion that $\k_2 = \c$. 
The second statement in (2) now follows from (2) of Corollary \ref{borel-tits-app-cor-1}.
\hfill$\Box$

We next study homomorphisms from real Lie groups and proceed towards proving Theorem \ref{real}.
We first need to establish an analogue of the {\it complete Jordan decomposition} for real Lie group $H$ which are covering groups of
$A(\r)^*$  for some algebraic group $A$ defined over $\r$.
We remark that Lemma \ref{jordan} is also established in \cite[Corollary 2.5.]{B1} with the additional restriction 
that $H$ is a finite cover of $A(\r)^*$.
In this set-up an analogous fact is proved in \cite[Example 4.11]{Ch2} on the usual Jordan decomposition of $H$. 

\begin{lemma}\label{jordan} 
Let $H$ be a real Lie group and $A$ be an algebraic group defined over $\r$.
Let $\eta : H \to A(\r)^*$ be a covering group homomorphism.
Let $h \in H$. Then there a unique pair of elements $s, w \in H$ with the following property
\begin{enumerate}
\item $h = s \, w = w \, s$.
\item $w \in \exp_H \circ {\rm d}\eta^{-1} ({\mathcal P}_{L(A(\r)^*)})$ and $\eta (s) = \eta(h)_e$
\end{enumerate}
\end{lemma}

\noindent{\bf Proof.} Let us denote ${\mathcal P}_H := \exp_H \circ {\rm d}\eta^{-1} ({\mathcal P}_{L(A(\r)^*)})$.
Using Lemma \ref{bijection-positive} we first conclude that $\eta ({\mathcal P}_H) = {\mathcal P}_{L(A(\r)^*)}$ and that
$\eta : {\mathcal P}_H \to {\mathcal P}_{L(A(\r)^*)}$ is a bijection. 
The uniqueness of the decomposition of $h$ is clear from this bijection. It now remains to show the existence.
Observe that there is unique $w \in {\mathcal P}_H$ such that
$\eta(w) = \eta(h)_p$. Clearly $\eta(h)_e = \eta ( s' )$ for some $s'  \in H$. Hence there is a $\gamma \in \ker \eta \subset Z(H)$
so that $ h = \gamma s' w$. We set $s := \gamma s'$. It is enough to show that $ s w s^{-1} = w$. 
It is easy to see that ${\mathcal P}_H$ is an invariant subset under conjugation action of $H$ on itself.
In view of the bijection $\eta : {\mathcal P}_H \to {\mathcal P}_{L(A(\r)^*)}$ we need to show that $\eta ( s w s^{-1}) = \eta(w)$.
It is now immediate that $\eta ( s w s^{-1}) =  \eta(h)_e \eta(h)_p \eta(h)_e^{-1} = \eta(h)_p= \eta(w)$.
\hfill$\Box$

We need the following well-known result.

\begin{thm}\label{center}{\rm (cf. \cite[Theorem 1.2, p 189]{Ho})} 
Let $B$ be a connected real Lie group. Then $Z(B)$ is contained in a closed, connected abelian subgroup of $B$.
\end{thm} 

An immediate known consequence of the above result is that $Z(B) \subset \exp_B (L(B))$.

\noindent{\bf Proof of Theorem \ref{real}.} 
We will divide the proof of the first part into three cases.

{\it Case 1:} We first prove the theorem in the case when the group $H$ is a covering group of $A(\r)^*$ for some algebraic group $A$ defined over $\r$. 

It is enough to prove that for any element $h \in H$ the element $\phi(h)$ is a unipotent element
in $G(\k)$. Let $\eta : H \to A(\r)^*$ be a covering homomorphism. 
In view of Lemma \ref{jordan} there exist a unique pair of elements $s, w \in H$ such that
$h = s \, w = w \, s$,
$w \in \exp_H \circ {\rm d}\eta^{-1} ( {\mathcal P}_{L(A(\r)^*)})$ and $\eta (s) = \eta(h)_e$.
As $\phi (h) = \phi (s) \, \phi (w)$ and as $\phi (s),  \, \phi (w)$ commute it is enough to show that
both  $\phi (s),  \, \phi (w)$ are unipotent elements in $G(\k)$.

As $w = \exp_H (X)$ for some $X \in {\rm d}\eta^{-1} ({\mathcal P}_{L(A(\r)^*)})$, it follows immediately that
$w$ is a divisible element in $H$. Hence $\phi (w) \in G(\k)$ is a divisible element. Now we use
Theorem \ref{unipotent} to conclude that $\phi (w) \in G(\k)$ is a unipotent element.

We next recall that $\eta (s) = \eta(h)_e$ lies in a maximal compact subgroup of $A(\r)^*$. As any maximal compact subgroup
of $A(\r)^*$ is connected there is a $Y \in L(H)$ such that $\eta (s) = \eta (\exp_H (Y))$. Hence $ s = \gamma \exp_H (Y)$
for some $\gamma \in \ker \eta \subset Z (H)$. Using Theorem \ref{center} it follows that $\gamma \in H$ is divisible and in particular
$\phi (\gamma) \in G(\k)$ is divisible. We now use Theorem \ref{unipotent} to conclude that 
$\phi (\gamma) \in G(\k)$ is a unipotent element. Again arguing as above we see that $\phi (\exp_H (Y)) \in       
G(\k)$ is a unipotent element. Clearly, $\phi (\exp_H (Y))$ and $\phi (\gamma)$ commute and in particular it follows that
$\phi (s) = \phi (\gamma) \, \phi (\exp_H (Y))$ is a unipotent element.  

{\it Case 2:} We now consider the case when $H$ is connected solvable group. Let $n := \dim H$.
We need the fact from  \cite[Theorem 3.18.11, p. 243]{V} that there exits Lie subgroups $H_1, \cdots , H_n$ and
one parameter Lie subgroups $C_1, \cdots ,C_n$ 
such that $H_1 = C_1$, $H_n = H$ and $H_i = C_1 \cdots C_i$, for all $1 \leq i \leq n$.
We may further assume that
$H_i$ is normal in $H_{i+1}$, for all $1 \leq i \leq n+1$. 
In particular $H_{i+1} = C_{i+1} H_i$, for all $1 \leq i \leq n+1$. 

We will prove that for all $i$ the group ${\rm Zcl}\, \phi (H_i)$ is a unipotent group. Clearly as $C_j$ is one parameter group, for all $j$, it is divisible and 
thus using Theorem \ref{unipotent} we conclude that ${\rm Zcl}\, \phi (C_j)$ is a unipotent subgroup. 
In particular ${\rm Zcl}\, \phi (H_1)$ is a unipotent subgroup.
Fixing a $i$ we now assume 
${\rm Zcl}\, \phi (H_i)$ is a unipotent group and will show that ${\rm Zcl}\, \phi (H_{i+1})$ is a unipotent group.
As $H_i$ is normal in $H_{i+1}$ it follows immediately that ${\rm Zcl}\,\phi (C_{i+1})$ normalises ${\rm Zcl}\, \phi (H_i)$.
Therefore, as ${\rm Zcl}\,\phi (C_{i+1}$ is unipotent subgroup, we see that  ${\rm Zcl}\,\phi (C_{i+1}) {\rm Zcl}\, \phi (H_i)$ is a unipotent subgroup.
Finally, as  ${\rm Zcl}\, \phi (H_{i+1}) \subset {\rm Zcl}\,\phi (C_{i+1}) {\rm Zcl}\, \phi (H_i)$, we conclude that
${\rm Zcl}\, \phi (H_{i+1})$ is a unipotent group.

{\it Case 3:} We now consider the general case when $H$ is a arbitrary connected Lie group. Let $R(H)$ be the radical of $H$.
Thus $R(H)$ is the maximal connected normal solvable subgroup of $H$. 
Then we have a (semisimple)Levi subgroup $M$ of $H$ such that $G = M \, R(H)$; see \cite[Theorem 3.18.13]{V}. As $M$ is semisimple
it is covering group of $A(\r)^*$ for some algebraic group defined over $\r$. So appealing to Case 1 we get that $\phi(M)$ lies in
a unipotent subgroup of $G$. Now the property $M= [M,M]$ implies that $\phi(M)=e$. So $\phi (H) = \phi (R(H))$. Now appealing to Case 2
we conclude that ${\rm Zcl}\,\phi (H) = {\rm Zcl}\,\phi (R(H))$ is a unipotent group in $G$.

This completes the proof of the first part. 

Now if $H'$ is an algebraic group defined over $\r$ then 
the semisimple Levi decomposition says that 
there is a semisimple algebraic subgroup $S$ over $\r$ such that $H' = S R(H')$
where $R(H)$ is the solvable radical of $H'$. Thus $H'(\r)^* = S (\r)^* R(H')(\r)^*$. As $S(\r)^*$ is a real semisimple Lie group
$\phi (S (\r)^*)=e$. Hence $\phi (H'(\r)^*) = \phi ( R(H')(\r)^*)$. Now as $R(H')$ is solvable algebraic group over $\r$ it follows
from $(2)$ of Theorem \ref{thm-solv} that $\phi ( R(H')(\r)^*)$ is abelian.
This completes the proof of the theorem
\hfill$\Box$

We now prove Theorem \ref{alg-closed}  on abstract homomorphisms from algebraic group over algebraically closed fields of arbitrary characteristic.

\noindent{\bf Proof of Theorem \ref{alg-closed}.} 
{\it Proof of (1)}: Let $B$ be a Borel subgroup of $H$. Then using (1) of Theorem \ref{thm-solv} it follows that
${\rm Zcl}\, \phi (B)$ is a unipotent subgroup of $G$.
Now let $h \in H$ be a arbitrary element. Then there is a $\alpha \in H$ such that $\alpha h \alpha^{-1} \in B$.
So from the above $\phi ( \alpha h \alpha^{-1})$ is a unipotent element in $G$, which in turn implies that
$\phi (h )$ is a unipotent element. This forces 
${\rm Zcl}\,\phi(H)$ to be a unipotent subgroup of $G$. 
Now as in the last part of proof of Theorem \ref{real} we see that $\phi (H) = \phi (R(H))$ where $R(H)$ is the solvable radical of 
$H$. Now we apply (1) of Theorem \ref{thm-solv} to see that $\phi (H)$ is abelian.

{\it Proof of (2):}
First note that this follows immediately from Proposition \ref{borel-tits-app}, which 
is a consequence of the Theorem \ref{borel-tits}. 
We next give another proof which is does not use
Proposition \ref{borel-tits-app}, and thus it is independent of  Theorem \ref{borel-tits}.
We divide the proof in two cases according to the characteristic of the field $\l$.

We consider the case when ${\rm Char}\l =0$. Using part (1) we see that   
${\rm Zcl}\,\phi(H)$ is a unipotent subgroup of $G$. But as $H$ is semisimple $H$ is perfect, that is, $[H,H]=H$.
As ${\rm Zcl}\,\phi(H)$ is nilpotent as a group it is now clear that $\phi$ is trivial.

We consider the case when ${\rm Char}\l >0$. Let ${\rm Char}\l = q$.
First observe that there is an integer $r$ such that for any unipotent element $u \in H$ one has
$u^{q^r} = e$.
Let $B$ be a Borel subgroup of $H$. As in the proof of (1) we will prove that ${\rm Zcl}\,\phi (B)$ is
a unipotent subgroup of $G$. As before, let $T$ be a maximal torus of $B$ and $U$ be the unipotent radical of $B$.
As $T$ is divisible, using Theorem \ref{unipotent} we conclude that
${\rm Zcl}\,\phi (T)$ is a unipotent subgroup of $G$.
On the other hand,  we have
$$
{\rm Zcl}\,  \phi (U) \subset \{ \,\, g \in G \,|\, g^{p^r} = e \,\,\}.
$$  

Hence all the elements of ${\rm Zcl}\,  \phi (U)$ are of order dividing $p^r$. Recall that in an algebraic group
over an algebraically closed field of characteristic zero
if all the elements are of finite order then it must be a finite group.
Thus ${\rm Zcl}\, \phi (U)$ is a finite group. 
Note that ${\rm Zcl}\,\phi (T)$ normalises ${\rm Zcl}\,  \phi (U)$. As ${\rm Zcl}\,\phi (T)$ is
a unipotent subgroup of $G$ and  ${\rm Char}\k =0$ it follows that 
${\rm Zcl}\,\phi (T)$ is a connected group. This implies that ${\rm Zcl}\,\phi (T)$ centralises ${\rm Zcl}\,  \phi (U)$. 
In particular for all $t\in T, u \in U$ the elements $\phi (t)$ and $\phi(u)$ commute.

Now consider the root system of $H$ with respect to the maximal torus $T$. Let $X(T)$ be the group of
characters of $T$ and let $\Delta \subset X(T)$ be the set of positive roots induced by the Borel subgroup $B$.
For $ \lambda \in \Delta$ let the associated root-group in $U$ be denoted by $U_\lambda$. 
It is straightforward, as also shown in the proof of \cite[Theorem A]{KM}, that if $x \in U_\lambda$
then there elements  $t, s \in T$ such that
\begin{equation}\label{alg-closed-1}
s x s^{-1} = x^{-1} t x t^{-1}.
\end{equation}
Since $\phi(x)$ and $\phi(t)$ commute, from the equation \eqref{alg-closed-1} we conclude that $\phi (x) =e$.
Now since the roots groups $U_\lambda$, for $ \lambda \in \Delta$, generate $U$ it follows that
$\phi (U) =e$. Thus ${\rm Zcl}\,\phi (B) = {\rm Zcl}\,\phi (T)$ is a unipotent group.
By arguing as in (1) we see that ${\rm Zcl}\,\phi (H)$ is a unipotent group. Now as $[H, H] =H$
it is clear that $\phi$ is trivial.
\hfill$\Box$

Let $l$ be a prime so that $l \neq p$.
We now consider abstract homomorphisms from the groups of rational points over $\Q_l$ to the
groups of rational points over $\q$.
We first deal with the groups of rational points of $\Q_l$-isotropic semisimple groups proving Theorem \ref{isotropic}.
We give two proofs; one using Corollary \ref{borel-tits-app-cor-3} and Lemma \ref{divalg}, and the other using
the well-known Jacobson-Morozov theorem (cf. Theorem \ref{jacobson-morozov}) and Proposition \ref{prop-solv} (3).

\noindent{\bf Proof of Theorem \ref{isotropic}.} 
We begin by proving $\phi ( H^+ (\Q_l)) = e$ which can be carried out in two different ways.

To see the first proof, using Lemma \ref{divalg}, we observe that, if $l \neq p$ then $\Q_l$ can not be embedded
in any finite extension of $\q$. In view of this the proof follows immediately from Corollary \ref{borel-tits-app-cor-3}.

To see the second proof we note that,
as  $H^+ (\Q_l)$ is generated by unipotent elements of $H(\Q_l)$ it is enough to prove 
$\phi (u) = e$ for any non-trivial unipotent element
$u \in H(\Q_l)$. Using Jacobson-Morozov theorem (see Theorem \ref{jacobson-morozov}) 
it follows that there is an abstract homomorphism
$\delta : {\rm SL}_2 ( \Q_l) \to H ( \Q_l) $ such that  
$$
\delta ( 
\begin{pmatrix}
1 & 1\\
0 & 1
\end{pmatrix}) = u.
$$

Let $B$ be the subgroup of ${\rm SL}_2 ( \Q_l)$ consisting of upper triangular matrices. 
Then using Proposition \ref{prop-solv} (3) we see that $\phi \circ \delta (B) = e$. In particular we have $\phi (u) =e$.  
Thus $\phi ( H^+ (\Q_l)) = e$.

{\it Proof of (1):}
This follows immediately from the fact that if $H$ is simply connected then $H(\Q_l) = H^+ (\Q_l)$; see Theorem \ref{kneser-tits} (3).

{\it Proof of (2):}
To prove (2) we just need to recall that when $H$ is a semisimple group
over $\Q_l$ then the quotient group $H(\Q_l)/ H^+ (\Q_l)$ is a finite abelian group; see Theorem \ref{kneser-tits} (2).
\hfill$\Box$

We next give  proof of Theorem \ref{sl(D)}. 
As ${\rm SL_1} (D)$ is
the group of rational points of a $\Q_l$-anisotropic algebraic group we need to employ results from \cite{Ch1} to 
handle this situation.                      

\noindent{\bf Proof of Theorem \ref{sl(D)}.} 
We start with certain facts regarding a specific filtration of ${\rm SL_1} (D)$  which asserts that
there are subgroups $\{U_i \}_{1=0}^{\infty}$ of ${\rm SL_1} (D)$ with $U_0 = {\rm SL_1} (D)$,  $U_{i+1} \subset U_i$, for all $i$ such that $U_i$ is a normal open subgroup of ${\rm SL_1} (D)$ for all $i$. This filtration was
investigated in \cite{Ri} but the relevant details can also be found in \cite[Section 4]{Ch1} or in \cite[Section 1.4, p. 27-33]{PR}.  We further recall the facts that $U_0 / U_1$ is a cyclic group and, for all $i \geq 1$, the group $U_i / U_{i+1}$ is isomorphic to the additive group of some finite dimensional vector spaces over certain finite extensions of the finite field 
${\mathbb Z}/l{\mathbb Z}$.  
Moreover, form the proof of \cite[Proposition 4.4]{Ch1} it follows that 
\begin{gather}\label{sl(D)-1}
| U_i / U_{i+1}|
= l^{rd}, \,\,\, {\rm if} \,\,\, i \not\equiv 0 \pmod{d} \nonumber\\
 \text{and}\,\,\, |U_i / U_{i+1}|
= l^{r(d-1)}, \,\,\, {\rm if} \,\,\, i \equiv 0 \pmod{d}.\nonumber
\end{gather}
In particular, the group $U_1$ is prosolvable as well as pro-$l$. 
Since $l \neq p$, using \cite[Theorem 3.3]{Ch1}, we conclude that $P_p : U_1 \to U_1$
is surjective. This is the same as saying that $U_1$ is $p$-divisible. We now appeal to Proposition \ref{p-divisible} to see that
${\rm Zcl}\,\phi (U_1) = F W$ for some finite group $F \subset G (\k)$ and a unipotent $\k$-subgroup $W$ of $G$ 
which is normalised by $F$. It further follows that
that $F$ is a (abstract) quotient group of the group $\phi (U_1)$ and hence a quotient of $U_1$. As $U_1$ is prosolvable, applying Theorem \ref{prosolv}, we see
that $F$ is a finite solvable group. Thus ${\rm Zcl}\,\phi (U_1) = F W$ is a solvable algebraic group. 
We define the derived series of $U_1$ by $C_1 := U_1$ and $C_j := [C_{j-1}, C_{j-1}]$ for all $j \geq 2$. 
As $U_1$ is normal in ${\rm SL_1} (D)$ so are the subgroups $C_j$, for all $j \geq 2$. Moreover, it
is clear that $C_j$ is non-central in ${\rm SL_1}(D)$, for all $j \geq 1$. Hence by Theorem \ref{riehm}
it follows that $C_j$ is open for all $j \geq 1$.
Clearly as ${\rm Zcl}\,\phi (U_1) $ is a solvable algebraic group there is an integer $r\geq 1$ such that $\phi (C_r) = (e)$. Now as $C_r \subset \ker \phi$ and
as $C_r$ is open the group  $\ker \phi$ is also an open subgroup. Thus the group the quotient $U_0/ \ker \phi$ is finite and
in particular we have that $\phi (U_0) $ is a finite group.
We now use Theorem \ref{prosolv} to see that $\phi (U_0) $ is a finite solvable group.
Recall that by \cite[Proposition 4.4]{Ch1} we have that 
\begin{equation}\label{sl(D)-2}
{\rm Ord} ( {\rm SL}_1 (D))= ( 1 + l^r + \cdots + (l^r)^{d-1}) l^\infty.
\end{equation}
Applying Theorem \ref{profinite} to the equality in \eqref{sl(D)-2} we conclude that the order of the group $\phi (U_0) $ divides 
the integer $( 1 + l^r + \cdots + (l^r)^{d-1}) l^\infty$. Thus 
$|\phi (U_0)| $ divides $( 1 + l^r + \cdots + (l^r)^{d-1}) l^m$, for some integer $m$.
This completes the proof.
\hfill$\Box$

\begin{rmk}
{\rm The methods of the proof of Theorem \ref{sl(D)} also apply to establish the following result.
Let $H$ be a profinite group and
$\phi : H \to G ( \k)$ be an abstract group homomorphism.
Assume that $H$
admits an open prosolvable subgroup $H_1$ such that ${\rm Ord} (H_1)$
is prime to $p$ and that the derived series of $H_1$, defined inductively by $C_1 := H_1$ 
and $C_j := [C_{j-1}, C_{j-1}]$ for all $j \geq 2$, consist of open subgroups.
Then $\phi (H)$ is a finite subgroup of $G ( \k)$.}
\end{rmk}

We now club Theorem  \ref{isotropic} and Theorem \ref{sl(D)} to deduce Theorem \ref{final-l-adic}.

\noindent{\bf Proof of Theorem \ref{final-l-adic}.}
{\it Proof of (1):} We first assume that $H$ is simply connected. Then $H$ is $\Q_l$-isomorphic to 
a direct product of finitely many $\Q_l$-simple simply connected groups $H_1, \cdots ,H_n$.
In particular $H (\Q_l)$ is isomorphic to the direct product $H_1 (\Q_l) \times \cdots
\times H_n (\Q_l)$. 
We next recall the well known fact that if $A$ is a simply connected, $\Q_l$-simple, $\Q_l$-anisotropic group then $A(\Q_l)$ is isomorphic to ${\rm SL}_1 (D)$ for some
central division algebra $D$ over a finite extension of $\Q_l$. 
In view of this fact, for the $\Q_l$-anisotropic factors in the above decomposition, we use 
Theorem \ref{sl(D)} and for the $\Q_l$-isotropic factors we use Theorem \ref{isotropic}
to conclude that $\phi ( H(\Q_l))$ is a finite solvable subgroup of $G (\k)$.

We now deal with the general case. Let $\tilde{H}$ be the simply connected covering group over $\Q_l$
and $\pi : \tilde{H} \to H$ be the covering homomorphism defined over $\Q_l$.  
Then 
using Corollary 2 of Section 5.6
and  Proposition 43 of Section 5.7 of \cite{Se2}
we have an exact sequence :
$$
1 \to {\rm Ker}\pi (\Q_l) \to \tilde{H} (\Q_l)  
\xrightarrow{\pi} H (\Q_l) \to H^1 ( \Q_l, {\rm Ker}\pi),
$$
where $H^1 (\Q_l, {\rm Ker}\pi)$ denotes the 
first Galois-cohomology group
of the Galois group ${\rm Gal}({\overline{\Q}}_l | \Q_l)$ with 
coefficients in the finite abelian group ${\rm Ker}\pi$. As $H^1 ( \Q_l, {\rm Ker}\pi)$ is finite abelian
so is the quotient group $H (\Q_l)/ \pi ( \tilde{H} (\Q_l))$. Now by the first case as above
$\phi \circ \pi ( \tilde{H} (\Q_l))$ is a finite solvable subgroup of $G (\k)$. This leads us to
conclude that $\phi ( H(\Q_l))$ is a finite solvable subgroup of $G (\k)$.

{\it Proof of (2):} The first statement of (2) follows from the (semisimple)-Levi decomposition of $H$.
Let $L$ be a semisimple Levi subgroup of $H$ defined over $\Q_l$. Then $H (\Q_l) = L (\Q_l) R(H) (\Q_l)$.
Now applying (1) we see that $\phi (L (\Q_l))$ is (finite)-solvable. This proves the first statement as
$R(H) (\Q_l)$ is a solvable group. 

To prove the second part of (2) we see that if $L$ is as above
then $H (\Q_l) = L (\Q_l) R_u ( H) (\Q_l)$. Now it follows from Corollary \ref{main-unip-cor} that
${\rm Zcl}\, \phi (R_u (H) (\Q_l))$ is abelian  and  unipotent. This completes the proof.
\hfill$\Box$

%\section*{Acknowledgments} 

\end{document}